\documentclass[12pt]{article}

\usepackage{amsmath,amssymb,amsfonts,amsthm}
\usepackage{graphics}
\usepackage{epsfig}
 \usepackage{amsmath}
\usepackage{amsfonts}
\usepackage{float}
\allowdisplaybreaks
\font\elevensf=cmss10 scaled\magstephalf

\newtheorem{theorem} {{\elevensf THEOREM}}[section]
\newtheorem{proposition} {{\elevensf PROPOSITION}}[section]

\newtheorem{lemma} {{\elevensf LEMMA}}[section]
\newtheorem{corollary} {{\elevensf COROLLARY}}[section]
\newtheorem{example} {{\elevensf EXAMPLE}}[section]

\newtheorem{remark} {{\elevensf REMARK}}[section]

\newcommand{\bpr}{\begin{proof}}
\newcommand{\epr}{\end{proof}}
\newcommand{\lb}{\left(}
\newcommand{\rb}{\right)}

\renewcommand\qed{$\blacksquare$}

\def\TagOnRight

\catcode`\@=11

\def\theequation{\@arabic{\c@section}.\@arabic{\c@equation}}

\begin{document}

\baselineskip 14pt
\parindent.4in
\catcode`\@=11

\begin{center}

{\Huge \bf Inverse Diffusivity Problem via Homogenization Theory} \\[5mm]

{\bf Tuhin GHOSH, Venkateswaran P.KRISHNAN and  Muthusamy VANNINATHAN }\\[4mm]

\textit{Centre for Applicable Matematics, Tata Institute of Fundamental Research, India.}\\[2mm]

Email : \textit{vanni@math.tifrbng.res.in\ ,\ vkrishnan@math.tifrbng.res.in\ ,\ tuhin@math.tifrbng.res.in }

\end{center} 

\begin{abstract}
\noindent
Polarization tensor corresponding to near zero volume inhomogeneities was
introduced in the pioneering work  by Capdeboscq-Vogelius \cite{CV1,CV2}. A
beautiful application of the polarization tensor to an inverse problem
involving inhomogeneities was also given by them. In this article, we take
an approach toward polarization tensor via homogenized tensor.
Accordingly, we introduce polarization tensor corresponding to
inhomogeneities with positive volume fraction.A relation between this
tensor and the homogenized tensor is found. Next, we proceed to examine the sense in which 
this tensor is continuous as the volume fraction tends to zero. 
Our approach has its own advantages, as we will see. In particular, it provides another method 
to deduce optimal estimates on polarization tensors in any dimension 
from those on homogenized tensors, along with the information on underlying microstructures.

\end{abstract}

\section{Introduction and statement of main results}

We consider a conducting object that occupies a bounded open set $\Omega \subset \mathbb{R}^N $ with smooth boundary. 
Let $\gamma_0$ denotes the constant background conductivity of the object in the absence of any inhomogeneities.
Let $\omega_\epsilon$ denote a set of inhomogeneities inside $\Omega$, we assume that the set $\omega_\epsilon$ is measurable for 
each $\epsilon$ and separated away from the boundary, that is, $\mbox{dist}(\omega_\epsilon , \partial \Omega) \geq d_0 > 0$. 
We assume that $|\omega_\epsilon| \to \delta \geq 0 $ as $\epsilon\to 0$. Such a situation was considered 
by Capdeboscq-Vogelius in  \cite{CV1,CV2} with $\delta =0$.
 
Let $\gamma_\epsilon$ denote the two-phase ($\gamma_0,\gamma_1$) conductivity profile of the medium in the presence of inhomogeneities,
that is
\begin{equation}\label{Conductivity_profile}
\gamma_ \epsilon(x)=\ \gamma_1\chi_{\omega_{\epsilon}}(x) + \gamma_0(1-\chi_{\omega_{\epsilon}}(x))\ \ x\in\Omega. 
\end{equation}
We assume that $ 0< \gamma_1 < \gamma_0 < \infty.$ ( The case $\gamma_0 < \gamma_1 $ can be treated analogously ).

The voltage potential in the above two phase medium is denoted by $u_{\epsilon}(x)$. It is the solution to the following two phase homogenization problem:
\begin{equation}\label{Inhomogeneities-BVP}
 \begin{aligned}
  \nabla\cdot(\gamma_{\epsilon}(x)\nabla u_{\epsilon}(x)) &=\  0 \mbox{ in } \Omega \\
     u_{\epsilon}(x) &=\ f(x) \textrm{ on } \partial \Omega\ \ \ \mbox{ for } f\in H^{1/2}(\partial \Omega).
 \end{aligned}
 \end{equation}
Let $M(\gamma_1,\gamma_0;\Omega)$ denote the set of all real symmetric positive definite matrices 
lying between $\gamma_1$ and $\gamma_0$. 
Assume that there exist $\theta\in L^\infty(\Omega;[0,1])$ and $\gamma^{*}\in M(\gamma_1,\gamma_0;\Omega)$, 
such that
\begin{equation}\label{ED7} \chi_{\omega_\epsilon}(x)\rightharpoonup\theta(x)\ \ \mbox{ weak* in }\ \ L^\infty(\Omega;[0,1])\end{equation}
and 
\begin{equation}\label{ED9}
\gamma_{\epsilon}(x)I\  \mbox{ H-converge to }  \ \gamma^{*}(x) \in M(\gamma_1,\gamma_0;\Omega) 
\end{equation}
in the sense that,
\begin{equation}\label{ED8}
\begin{aligned}
 u_{\epsilon} &\rightharpoonup u \quad\mbox{weakly in }H^{1}(\Omega),\\
 \gamma_{\epsilon}\nabla u_{\epsilon} &\rightharpoonup \gamma^{*}\nabla u \quad\mbox{weakly in }(L^{2}(\Omega))^N
\end{aligned}
\end{equation}
where $u$ is the solution of the homogenized equation : 
\begin{equation}\label{Homogeneities-BVP}
\begin{aligned}
 \nabla\cdot(\gamma^{*}\nabla u(x))&=\ 0 \quad\mbox{in }\Omega \\
 u(x)&=\ f(x)  \quad\mbox{on }\partial\Omega.
\end{aligned}
\end{equation}
It is known above convergences \eqref{ED7},\eqref{ED9} hold for a suitable subsequence. Our hypothesis here is that they hold for the entire sequence.
Motivated by the inverse problem of determining the measure of dilute inhomogeneities, the authors of 
\cite{CV1,CV2} consider the equation in \eqref{Inhomogeneities-BVP} with Neumann boundary condition and 
the asymptotic observation/measurement which is nothing but 
the difference in potential on the boundary. Since we have imposed the Dirichlet boundary condition
in the problem \eqref{Inhomogeneities-BVP}, we consider the following
asymptotic observation/measurement which is nothing but the current perturbation on the boundary of the domain :
\begin{equation}\label{Boundary current perturbation}
 (\gamma_{\epsilon}\frac{\partial u_{\epsilon}}{\partial \nu}- \gamma_0\frac{\partial u}{\partial \nu}) |_{\partial \Omega} \mbox{ as } \epsilon \rightarrow 0,
 \end{equation}
 where $\nu$ is the outer unit normal to $\partial \Omega$. 
 A general asymptotic formula for boundary current perturbations under the condition that the volume of inhomogeneities goes 
 to zero was derived in \cite{CV1}. 
 For earlier works in this direction under certain restrictions on the inhomogeneities, 
 we refer \cite{VV,BMV,AMV,BEV}. 
 The asymptotic formula derived in \cite{CV1,CV2} involved the so-called polarization tensors which form a set of 
  macro coefficients associated to dilute inhomogeneities. They 
 studied several properties of these tensors in \cite{CV1,CV2,CV2A,CV3,CV4}. 
 For a comprehensive treatment of polarization tensors, we refer the reader 
 to \cite{AK}.\\
 
 In this article, as already mentioned, we consider inhomogeneities with total volume $\delta > 0$. The 
 appropriate quantity for which we seek the asymptotic formula is then the following :
\begin{equation}\label{Boundary current perturbation-2}
 (\gamma_{\epsilon}\frac{\partial u_{\epsilon}}{\partial \nu}- \gamma^{*}\frac{\partial u}{\partial \nu}) |_{\partial \Omega} \mbox{ as } \epsilon \rightarrow 0,
 \end{equation}
This generalizes \eqref{Boundary current perturbation} because when $\delta=0$ the homogenized tensor 
$\gamma^{*}$ coincides with $\gamma_0I$. Naturally the above formula involves a new polarization tensor corresponding to 
$\delta $ positive.
Roughly speaking, polarization tensor provides first order 
approximation of the homogenized tensor, $\gamma^{*}$ for small
volume proportion of inhomogeneities \cite{AK,CV1,CMV,LR}.
The papers \cite{CV1,CV2}
present a microscopic interpretation of the polarization
tensor. To establish its regularity properties
near zero-volume fraction, we find it appropriate to introduce the polarization
tensor at non-zero volume fraction and then study its asymptotic properties as volume
fraction goes to zero. Our approach examines the precise sense in which polarization
tensor is continuous at zero volume fraction. It has also other advantages
as well. One advantage is that we obtain a relation between polarization tensor 
and homogenized tensor for non-zero volume fraction and another advantage is that it 
allows to use the knowledge about the homogenized tensor to deduce the properties of 
polarization tensor. As an example, we deduce  
optimal estimates on polarization tensors in any dimension from
those on homogenized tensors, along with the underlying microstructures.
At the core of our approach lies the following approximation result : any polarization 
tensor corresponding to zero volume fraction can be obtained as an appropriate limit of  
polarization tensor corresponding to non-zero volume fraction ( Theorem \ref{approximation} ).

We close this section by stating how the article is organized. In the reminder of this section we 
state our main results. The first theorem gives asymptotic expression for the perturbed current 
\eqref{Boundary current perturbation-2} as $|\omega_\epsilon|\rightarrow \delta$.
The proof of this Theorem is presented in Section 2. This naturally calls the introduction of the polarization
tensor denoted $M^\theta$ of the non-zero volume fraction. An important relation linking this
polarization tensor with homogenized tensor $\gamma^{*}$ is stated in Proposition \ref{PT and HT Proposition}.
Bounds on $M^{\theta}$ are easily deduced are those of $\gamma^{*}$. Theorem \ref{approximation} is concerned with 
a continuity property of polarization tensor whose proof is presented in Section 3. 
Theorem \ref{Optimal trace bounds theorem} presents the optimal bounds in terms of trace inequalities on polarization tensor
near zero volume fraction as a consequence of Theorem \ref{approximation}. Regarding the microstructures underlying polarization tensors $M^{0}$ with near zero volume 
inhomogeneities, we recall that  \cite{CV2A} shows that equality in the 
above trace bounds holds for the so-called ``washers'' microstructures
in two dimension. Numerical evidence for the same is provided in \cite{ACEK}. 
Examples of ``thin'' inhomogeneities are treated in \cite{CV3}.
Thanks to our approach, we are able to compute $M^{0}$ corresponding to sequential laminates of any rank  in any dimension 
with relative ease. Such tensors ``fill up'' the region in the phase space defined by the trace inequalities \eqref{LU-bounds}.
See Section \ref{Polarization tensors zero volume fraction optimal trace bounds}. Such a computation seems hard without passing through the homogenized tensor.\\

In order to state our result, we require a few preliminaries from homogenization theory \cite{A,T}.
The homogenized tensor is obtained from oscillating test functions which are defined by 
\begin{equation}\label{w-epsilon}
\begin{aligned}
-\nabla\cdot (\gamma_{\epsilon}(x)\nabla w_{\epsilon}^i(x)) &=\ -\mbox{div}(\gamma^{*}e_i) \quad\mbox{in }\Omega. \\
  w_{\epsilon}^i &=\ x_i \quad \mbox{ on }\partial\Omega.
\end{aligned}
\end{equation}
The matrix $W_{\epsilon}$ defined by its columns ($\nabla w_{\epsilon}^i)_{1\leq i \leq N}$ is called the corrector matrix with the following properties: 
\begin{equation}\label{Corrector Theory}
 \begin{aligned}
  W_{\epsilon} &\rightharpoonup I \textrm{ weakly in }  L^{2}(\Omega)^{N\times N} \\
 \gamma_{\epsilon}(x) W_{\epsilon}(x) &\rightharpoonup \gamma^{*} \mbox{ weakly in } L^{2}(\Omega)^{N\times N}.
 \end{aligned}
\end{equation}
We can write 
\begin{equation}\label{Corrector theory result}
\nabla u_{\epsilon} = W_{\epsilon}\nabla u + r_{\epsilon} 
\end{equation}
where $r_{\epsilon} \rightarrow 0 $ strongly in $(L^1_{\mbox{loc}}(\Omega))^N$.\\

We also need the so-called boundary Green's function which is defined as follows.  
For $y\in \partial\Omega$, consider 
\begin{equation}\label{Background-BVP}
 \begin{aligned}
  \nabla\cdot(\gamma_0\nabla_x D(x,y)) &=\ 0 \mbox{ in } \Omega \\
     D(x,y) &=\ \delta_y(x)  \mbox{ on }\partial \Omega. 
 \end{aligned}
\end{equation}
\begin{theorem}\label{Polarization tensor thm-a}
Let $\delta > 0$. Given $f\in H^{1/2}(\partial \Omega)$, let $u$ and $u_{\epsilon}$ denote the solutions to \eqref{Homogeneities-BVP} and  \eqref{Inhomogeneities-BVP} respectively. 
Then there exists a subsequence (still denoted by $\epsilon$), a regular positive compactly supported Radon measure $\mu^{\theta}$, and a matrix-valued function $M^{\theta} \in  L^2(\Omega, d \mu^\theta)$ ( called polarization tensor ) 
such that 
\begin{equation}
\begin{aligned}
(\gamma_{\epsilon}\frac{\partial u_{\epsilon}}{\partial \nu} - \gamma^{*}\frac{\partial u}{\partial \nu})(y)
&=\ |\omega_{\epsilon}|\int_{\Omega} (\gamma_1-\gamma_0) M^{\theta}_{ij}(x)\frac{\partial u}{\partial x_i}(x)\frac{\partial D}{\partial x_j}(x,y)d \mu^{\theta}(x) \\
&\quad + \int_{\Omega} (\gamma_0I - \gamma^{*}(x))_{ij}\frac{\partial u}{\partial x_i}(x)\frac{\partial D}{\partial x_j}(x,y)d x 
+ \mathrm{o}(1). 
\end{aligned}
\end{equation}
The $\mathrm{o}(1)$ term goes to zero uniformly in $y$ as $\epsilon$ goes to zero. 
\end{theorem}
\begin{remark}
 The polarization tensor $M^{\theta}$ depends on the microstructures $\omega_\epsilon$ under consideration.
 Considering $\gamma_{\epsilon}$ in the case of $\delta=0$ the above Theorem \ref{Polarization tensor thm-a} 
 was proved in \cite{CV1}.
\end{remark}
\begin{remark}\label{per-cons}
 If we consider periodic microstructures with a given volume fraction of inhomogeneities,
 then $\theta>0$ is constant which is equal to the volume fraction and the corresponding 
 polarization tensor $M^\theta$ is also constant. Moreover, the corresponding measure 
 $\mu^\theta$ is given by $\frac{1}{|\Omega|}dx$ ( see \eqref{ED5} ).
\end{remark}
\begin{theorem}\label{approximation}
For a given $d\mu^{0}$, $M^{0}$, the polarization tensor of near zero-volume fraction with $M^0\in L^2(\Omega,d \mu^{0})$, 
for any point $x_0\in support\ of\ d\mu^{0}$, $\mu^{0}$ almost everywhere  there exists a sequence 
$\theta_{x_0}^n \in (0,1] $ depending upon the point $x_0$ and a sequence of polarization tensors 
$M^{\theta_{x_0}^n}$ which are constant, such that as $n\rightarrow \infty$, $\theta_{x_0}^n \rightarrow 0$ and $ M^{\theta_{x_0}^n} \rightarrow M^{0}(x_0)$, $\mu^{0}$ almost everywhere $x_0$.
\end{theorem}
\begin{theorem}\label{Optimal trace bounds theorem}
\textbf{(a): }\ Let $M^{0}$ denote a polarization tensor 
corresponding to zero volume fraction. Then for $\mu^{0}$ almost everywhere $x\in support\ of\ \mu^{0}$,
we have the following pointwise trace bounds, 
\begin{equation}\begin{aligned}\label{LU-bounds}
&\mbox{Lower Bound : }\quad  \mathrm{trace}(M^0(x))^{-1} \leq\ (N-1) + \frac{\gamma_1}{\gamma_0}\\
&\mbox{Upper Bound : }\quad \mathrm{trace}(M^0(x)) \ \ \leq\ (N-1) + \frac{\gamma_0}{\gamma_1}.
\end{aligned}\end{equation}
\textbf{(b): }\ 
These bounds are optimal in the sense that any $(\lambda_1(x),..,\lambda_N(x))$ satisfying pointwise
\begin{align*}
\sum_{i=1}^N \lambda_i(x) &\leq\ (N-1) + \frac{\gamma_0}{\gamma_1}.\\
\sum_{i=1}^N \frac{1}{\lambda_i(x)}  &\leq\ (N-1) + \frac{\gamma_1}{\gamma_0}.
\end{align*}
arises as the eigenvalues of a polarization tensor of zero volume fraction at that point $x$.  
\end{theorem}

\section{Proof of Theorem \ref{Polarization tensor thm-a}}
\setcounter{equation}{0}
In this section, we obtain an asymptotic formula for the boundary current perturbations when the volume of inhomogeneities goes to $\delta>0$ as $\epsilon\to 0$.

We begin by applying the divergence formula a few times to obtain an expression for the boundary current perturbation.

Multiply $D(x,y)$ in \eqref{Inhomogeneities-BVP} and \eqref{Homogeneities-BVP} and using the divergence formula, we have 

\begin{align*}
   \int_{\Omega}\gamma_{\epsilon}(x)\nabla u_{\epsilon}(x)\cdot\nabla_x D(x,y)d x = \int_{\partial \Omega}\gamma_{\epsilon}(x) \frac{\partial u_{\epsilon}}{\partial \nu}(x)D(x,y) d \sigma(x)\\
  \int_{\Omega}\gamma^{*}(x)\nabla u(x)\cdot\nabla_x D(x,y)d x = \int_{\partial \Omega}\gamma^{*}(x) \frac{\partial u}{\partial \nu}(x)D(x,y) d \sigma(x)
 \end{align*}
Subtracting these two equations, we get
\begin{equation}\label{Rep-formula}
 \begin{aligned}
   &((\gamma_{\epsilon}(x)\frac{\partial u_{\epsilon}}{\partial \nu}(x) - \gamma^{*}(x)\frac{\partial u}{\partial \nu}(x)))D(x,y) |_{\partial \Omega} = \gamma_0(\frac{\partial u_{\epsilon}}{\partial \nu} - \frac{\partial u}{\partial \nu})(y)\\
&= \int_{\Omega}\gamma_{\epsilon}(x)\nabla u_{\epsilon}(x)\cdot\nabla_x D(x,y)d x - \int_{\Omega}\gamma^{*}(x)\nabla u(x)\cdot\nabla_x D(x,y)d x\\
&= \int_{\omega_{\epsilon}}(\gamma_1 - \gamma_0)\nabla u_{\epsilon}(x)\cdot\nabla_x D(x,y)d x + \int_{\Omega}\left[\gamma_0\nabla u_{\epsilon}-\gamma^{*}(x)\nabla u(x)\right]\cdot\nabla_x D(x,y)d x.
\end{aligned}
\end{equation}
Now multiplying \eqref{Background-BVP} by $u_{\epsilon}$ and $u$ and again applying the divergence formula, we have,
  \begin{align*}
&\int_{\Omega}\gamma_0\nabla_x D(x,y)\cdot\nabla u_{\epsilon}(x)d x = \int_{\partial \Omega}\gamma_0 \frac{\partial D(x,y)}{\partial \nu}(x)u_{\epsilon}(x) d\sigma(x)\\
&\int_{\Omega}\gamma_0\nabla_x D(x,y)\cdot\nabla u(x)dx = \int_{\partial \Omega}\gamma_0 \frac{\partial D(x,y)}{\partial \nu}(x)u(x) d\sigma(x).
\end{align*}
Since $u_{\epsilon}=u=f$ on $\partial \Omega$, we have
\[
\int_{\Omega}\gamma_0\nabla_x D(x,y)\cdot\nabla u_{\epsilon}(x)d x =\ \int_{\Omega}\gamma_0\nabla_x D(x,y)\cdot\nabla u(x)d x. 
\]
Therefore 
\begin{equation}\label{Current perturbation formula}
\begin{aligned}
(\gamma_{\epsilon}\frac{\partial u_{\epsilon}}{\partial \nu} - \gamma^{*}\frac{\partial u}{\partial \nu})(y)=\ \gamma_0(\frac{\partial u_{\epsilon}}{\partial \nu} - \frac{\partial u}{\partial \nu})(y)
&=\ \int_{\omega_{\epsilon}}(\gamma_1 - \gamma_0)\nabla u_{\epsilon}(x)\cdot\nabla_x D(x,y)\ dx \\
&\ \ +\int\limits_{\Omega} \lb \gamma_0I-\gamma^{*}(x)\rb \nabla u(x)\cdot \nabla_{x}D(x,y)\ dx.
\end{aligned}
\end{equation}
By \eqref{Corrector theory result}, we can rewrite  
\begin{equation}\label{First term of current perturbation formula}\begin{aligned}
\int_{\omega_{\epsilon}}(\gamma_1 - \gamma_0)\nabla u_{\epsilon}(x)\cdot\nabla_x D(x,y)d x &=\ \int\limits_{\omega_{\epsilon}} (\gamma_1 - \gamma_0)(W_{\epsilon}\nabla u+r_{\epsilon})\cdot \nabla_{x} D(x,y) d x\\
&=\ |\omega_{\epsilon}|\int\limits_{\Omega} (\gamma_1-\gamma_0)\frac{\chi_{\omega_{\epsilon}}(x)}{|\omega_{\epsilon}|} \frac{\partial w_{\epsilon}^{i}}{\partial x_{j}}\frac{\partial u}{\partial x_{j}}\frac{\partial D(x,y)}{\partial x_{i}}  d x\\
&\quad +\int\limits_{\omega_{\epsilon}}(\gamma_1-\gamma_0)r_{\epsilon}\cdot \nabla_{x} D(x,y) d x.
\end{aligned}\end{equation}
We first focus on \eqref{First term of current perturbation formula} above. We have the following results.
\begin{proposition} Let $\{ w_{\epsilon}\}$ be as in \eqref{w-epsilon}. Then 
\begin{equation}\label{ED1} \lVert\nabla w_{\epsilon}^{i} - e_i \rVert_{2}\ \leq\ C|\omega_{\epsilon}|^{\frac{1}{2}} + \lVert(\gamma^{*}-\gamma_0I)e_i\rVert_{2}.\end{equation}
\end{proposition}
\bpr 
Consider the corrector equation \eqref{w-epsilon} and multiply by $( w_{\epsilon}^{i} - x_i)$ and using divergence formula, we have
\[
\int_{\Omega}  \gamma_{\epsilon}\nabla  w_{\epsilon}^{i}\cdot (\nabla w_{\epsilon}^{i} -e_{i})\ d x=\ \int_{\Omega} \gamma^{*}(x)(\nabla w_{\epsilon}^{i}-e_{i})\ d x.
\]
Rewriting it as, 
\begin{align*}
\int\limits_{\Omega} \gamma_{\epsilon}(\nabla  w_{\epsilon}^{i} - e_i)\cdot(\nabla  w_{\epsilon}^{i} - e_i)\ d x 
\ =\ \int\limits_{\Omega} (\gamma^{*}-\gamma_0I)e_i\cdot(\nabla w_{\epsilon}^{i} - e_{i})\ d x&\\
- \int\limits_{\Omega}(\gamma_1 - \gamma_0)& \chi_{\omega_{\epsilon}}  e_{i}\cdot(\nabla w^{i}_{\epsilon} - e_i)\ d x.
\end{align*}
\hfill\epr
\begin{remark}\label{qq8}
 If $\delta = 0 $ then as $\gamma^{*} = \gamma_0I$, in that case \eqref{ED1} becomes $\lVert\nabla w_{\epsilon}^{i} - e_i \rVert_{2}\ \leq\ C|\omega_{\epsilon}|^{\frac{1}{2}}.$
 \end{remark}
\begin{lemma}\label{Matrix Lemma}
Let $\{w_{\epsilon}\}$ be as in \eqref{w-epsilon}. There exists a subsequence still denoted by $\epsilon$, function $\theta\in L^{\infty}(\Omega,[0,1])$, compactly supported a positive Radon measure $\mu^{\theta}$, and a 2-tensor $M\in L^{2}(\Omega, d \mu^{\theta})$ such that 
\[
\frac{1}{|\omega_{\epsilon}|}\chi_{\omega_{\epsilon}}(x) d x \stackrel{*}{\rightharpoonup} d \mu^{\theta} \mbox{ in }  (C^0(\overline{\Omega}))^{*},
\]
\[
\frac{1}{|\omega_{\epsilon}|}\chi_{\omega_{\epsilon}}(x) \frac{\partial w_{\epsilon}^i}{\partial x_j}d x \stackrel{*}{\rightharpoonup} M^{\theta}_{ij}d \mu^{\theta} \mbox{ in }  (C^0(\overline{\Omega}))^{*}.
\]
\end{lemma}
\bpr
 Since $\chi_{\omega_{\epsilon}}(x) \in  L^{\infty}(\Omega,\{0,1\})$ there exists a $\theta \in   L^{\infty}(\Omega,[0,1])$ and a subsequence such that
\[
\chi_{\omega_{\epsilon}}(x)\stackrel{*}{\rightharpoonup}  \theta(x) \mbox{ in }  L^{\infty}(\Omega,[0,1]).
\]
In particular, we have 
$$|\omega_{\epsilon}| = \int_{\Omega} \chi_{\omega_{\epsilon}}(x) d x \rightarrow \int_{\Omega} \theta(x) d x\ = \delta.$$
Now  since $\frac{1}{|\omega_{\epsilon}|}\chi_{\omega_{\epsilon}}(x)$ is bounded in $L^{1}(\Omega)$, from Banach-Alaoglu theorem and Riesz Representation theorem, there exists a regular, Radon measure $\mu^{\theta}$ and a subsequence (denoted by $\epsilon$) such that 
\begin{equation*}
 \frac{1}{|\omega_{\epsilon}|}\chi_{\omega_{\epsilon}}(x) \stackrel{*}{\rightharpoonup} d \mu^{\theta} \mbox{ in }(C^0(\overline{\Omega}))^{*} .
 \end{equation*}
And if $\delta > 0$ we see 
 \begin{equation}\label{ED5}
d \mu^{\theta} = \frac{\theta(x)}{\int_{\Omega} \theta(x)d x}d x \quad \mbox{almost everywhere}\ x \ \mbox{ in }\Omega.
\end{equation}

Now, $\frac{1}{|\omega_{\epsilon}|}\chi_{\omega_{\epsilon}}(x) \nabla w_{\epsilon}^i $ is also bounded $L^{1}(\Omega)$ for each $i=1,..,N$ as
\begin{align*} 
|| \frac{1}{|\omega_{\epsilon}|}\chi_{\omega_{\epsilon}}(x) \nabla w_{\epsilon}^i ||_{L^1(\Omega)} &= \int_{\Omega}  \frac{1}{|\omega_{\epsilon}|}\chi_{\omega_{\epsilon}}(x)|(\nabla w_{\epsilon}^i - e_i) + e_i|\ dx\\
&\ \leq \frac{1}{|\omega_{\epsilon}|^{\frac{1}{2}}} (\int_{\Omega} |\nabla w_{\epsilon}^i - e_i |^2\ dx)^{1/2} +  \frac{1}{|\omega_{\epsilon}|^{\frac{1}{2}}}( \int_{\omega_{\epsilon}} 1\ dx)^{1/2}\\
&\ \leq C  \quad\mbox{ by using }\eqref{ED1}\mbox{ for both }\delta=0 \mbox{ and }\delta > 0 \mbox{ cases. }   
\end{align*}
So there exists a subsequence and a Radon measure $d M^{\theta}_{ij}$ such that
\begin{equation}\label{ED2}
 \frac{1}{|\omega_{\epsilon}|}\chi_{\omega_{\epsilon}}(x) \frac{\partial w_{\epsilon}^i}{\partial x_j}d x \stackrel{*}{\rightharpoonup} d M^{\theta}_{ij} \mbox{ in }  (C^0(\overline{\Omega}))^{*}.
\end{equation}
Now, 
      \begin{align*}
       \lvert\int_\Omega \phi(x) d M^{\theta}_{ij} \rvert &= \lvert\hspace{2pt} \textrm{lim} \hspace{2pt} \int_{\Omega}  \frac{1}{|\omega_{\epsilon}|}\chi_{\omega_{\epsilon}}(x) \frac{\partial w_{\epsilon}^i}{\partial x_j} \phi(x) d x\rvert\\
&\leq \underbar{lim} \hspace{1pt} \left(\frac{1}{|\omega_{\epsilon}|}\right)^{1/2} \left( \int \left| \frac{\partial w_{\epsilon}^i}{\partial x_j}\right|^2 d x \right)^{1/2} \left( \int \frac{1}{|\omega_{\epsilon}|}\chi_{\omega_{\epsilon}}(x)|\phi|^2 d x\right)^{1/2}\\
& \leq C\left(\int \theta (x)d x \right)^{-1/2} \left(\int |\phi|^2 d \mu^{\theta}\right)^{1/2}  \mbox{ for all }  \phi \in C^0(\bar{\Omega}).\\
\end{align*}
Therefore, $\phi \longmapsto \int_{\Omega} \phi d M^{\theta}_{ij}$ can be extended to a bounded linear functional on $L^2(\Omega, d \mu^{\theta})$. Now  by Riesz-Representation theorem
$$ \int_{\Omega} \phi\ d M^{\theta}_{ij} =  \int_{\Omega} \phi\ M^{\theta}_{ij} d \mu^{\theta} \ \ \textrm{for some function}\ \ M^{\theta}_{ij} \in  L^2(\Omega, d \mu^{\theta}).$$
Hence, 
$$ d M^{\theta}_{ij}=M^{\theta}_{ij}d \mu^{\theta}.$$
For $\delta > 0$ we also see that
\begin{equation}\label{ED3} d M^{\theta}_{ij}=M^{\theta}_{ij}d \mu^{\theta} =\frac{M^{\theta}_{ij}(x)\theta(x)}{\int \theta(x)d x}\ d x \quad \mbox{ almost everywhere in }\Omega.\end{equation}
This completes the proof of the lemma.
\hfill\epr
\begin{remark}
We call the matrix $M^{\theta}= (M^{\theta}_{ij})$ polarization tensor with non-zero volume fraction.
For $\delta=0,$ we have $\gamma^{*}=\gamma_0I$, and the measure $\mu^0$ and the matrix $M^0$ are given by  ( cf. \cite{CV1} ).
\begin{equation}\label{m0-d}
\frac{1}{|\omega_{\epsilon}|}\chi_{\omega_{\epsilon}}(x) \frac{\partial w_{\epsilon}^i}{\partial x_j}dx \stackrel{*}\rightharpoonup d M^{0}_{ij}= M^{0}_{ij}d \mu^{0}  \textrm{ in } (C^0(\overline{\Omega}))^{*}
\end{equation}
with 
\begin{equation}\label{mu0-d}
\frac{1}{|\omega_{\epsilon}|}\chi_{\omega_{\epsilon}}(x) \stackrel{*}\rightharpoonup d \mu^{0} \mbox{ in } (C^0(\overline{\Omega}))^{*}\ \mbox{ and }\ M^0 \in L^2(\Omega,d \mu^{0}).
\end{equation}
The correctors $(w_{\epsilon}^i(x))_{1\leq i\leq N}$ are the solution of
\begin{equation}\label{ED6} -\nabla\cdot(\gamma_{\epsilon}(x)\nabla w_{\epsilon}^i(x))\ =\ - \mathrm{div}(\gamma_0e_i)\ \mbox{ in }\ \Omega, \quad w_{\epsilon}^i(x)\ =\ x_i\ \mbox{ on }\ \partial\Omega. \end{equation}
\end{remark}
\begin{remark}
For $\delta>0$, \eqref{ED2} shows that product of two weakly convergent sequence $\frac{1}{|\omega_{\epsilon}|}\chi_{\omega_\epsilon}dx$ and $\nabla w_{\epsilon}^i $ converges weakly
and the limit is not the product of limits. The measure $\mu^{\theta}$ is absolutely
continuous with respect to Lebesgue measure. $M^\theta$ is defined
almost everywhere in $\Omega$ with respect to the Lebesgue measure.
Such properties do not hold for $\mu^{0}, M^{0}$.
For $\delta =0 $, $M^0$ is defined over set of the support of $d\mu^{0}$. As a convention we define $M^0$ equal to identity elsewhere in $\Omega$. 
\end{remark}

We now prove Theorem \ref{Polarization tensor thm-a}.
\bpr[Proof of Theorem \ref{Polarization tensor thm-a}] Let $y\in \partial \Omega$. We recall the following two equations: 
\begin{equation}\label{Current perturbation formula-b}
\begin{aligned}
\gamma_0(\frac{\partial u_{\epsilon}}{\partial \nu} - \frac{\partial u}{\partial \nu})(y)
=&\underbrace{\int_{\omega_{\epsilon}}(\gamma_1 - \gamma_0)\nabla u_{\epsilon}(x)\cdot\nabla_x D(x,y)d x}_{J} \\
&+\int\limits_{\Omega} \lb \gamma_0I-\gamma^{*}(x)\rb \nabla u(x)\cdot \nabla_{x}D(x,y)\ d x
\end{aligned}
\end{equation}
\begin{align*}\label{First term of current perturbation formula-b}
J =&\int\limits_{\omega_{\epsilon}} (\gamma_1 - \gamma_0)(W_{\epsilon}\nabla u+r_{\epsilon})\cdot \nabla_{x} D(x,y) d x\\
=&\ |\omega_{\epsilon}|\int\limits_{\Omega} (\gamma_1-\gamma_0)\frac{\chi_{\omega_{\epsilon}}(x)}{|\omega_{\epsilon}|} \frac{\partial w_{\epsilon}^{i}}{\partial x_{j}}\frac{\partial u}{\partial x_{j}}\frac{\partial D(x,y)}{\partial x_{i}}\  d x
+\int\limits_{\omega_{\epsilon}}(\gamma_1-\gamma_0)r_{\epsilon}\cdot \nabla_{x} D(x,y)\ d x.
\end{align*}

As in \cite[pp. 169]{CV1}, we have the following: Let $\Omega_{1}\subset\subset \Omega$ denote a compact set that strictly contains the inhomogeneties $\omega_{\epsilon}$.  Given $y\in \partial \Omega $, we can find a vector-valued function $\phi_y \in C^{0}(\overline{\Omega}_1)$ such that $\phi_y(x) = \nabla_xD(x,y)$ for all  $x \in \Omega_1$. Also since $u$ is smooth in the interior of $\Omega$ and $d M_{ij}^{\theta}$ is supported in a compact subset of $\Omega$, we have
\[
\lim_{\epsilon\to 0} \int\limits_{\Omega} (\gamma_1-\gamma_0)\frac{\chi_{\omega_{\epsilon}}(x)}{|\omega_{\epsilon}|} \frac{\partial w_{\epsilon}^{i}}{\partial x_{j}}\frac{\partial u}{\partial x_{j}}\frac{\partial \phi_{y}(x)}{\partial x_{i}}  d x = \int\limits_{\Omega}(\gamma_1-\gamma_0) M_{ij}^{\theta}\frac{\partial u}{\partial x_{j}}\frac{\partial D(x,y)}{\partial x_{i}} d \mu^{\theta}.
\]
Since $r_{\epsilon} \to 0$ strongly in $(L^1_{\mathrm{loc}}(\Omega))^{N})$, we have 
\[
\int\limits_{\omega_{\epsilon}}(\gamma_1-\gamma_0)r_{\epsilon}\cdot \nabla_{x} D(x,y) d x \to 0 \mbox{ as } \epsilon\to 0.
\]
Hence 
\begin{align*}
\gamma_0(\frac{\partial u_{\epsilon}}{\partial \nu} - \frac{\partial u}{\partial \nu})(y)&=\ |\omega_{\epsilon}|\int\limits_{\Omega}(\gamma_1-\gamma_0) M_{ij}^{\theta}\frac{\partial u}{\partial x_{j}}\frac{\partial D(x,y)}{\partial x_{i}}\ d \mu^{\theta}\\
&\quad +\int\limits_{\Omega} \lb \gamma_0I-\gamma^{*}(x)\rb \nabla u(x)\cdot \nabla_{x}D(x,y)\ d x+\mathrm{o}(1)
\end{align*}
where the $o(1)$ goes to zero uniformly in $y$ as $\epsilon$ goes to zero. This completes the proof of Theorem \ref{Polarization tensor thm-a}.
\hfill\epr
\subsection{Some properties of polarization tensor}\label{Polarization tensor Properties Section}
\begin{proposition}[Relation between polarization and homogenization tensor]\label{PT and HT Proposition}
Let $\delta >0$. The polarization tensor and homogenization tensor are related as follows:
\begin{equation}\label{PT and HT relation}
\theta(x)(\gamma_1-\gamma_0)M^{\theta}(x)  = \gamma^{*}(x) - \gamma_0I.
\end{equation}
\begin{remark}
If $\delta=0$, then above equality holds trivially because $\theta=0$ almost everywhere and $\gamma^{*}=\gamma_0I$. In some cases, above relation then degenerates.
\end{remark}
\end{proposition}
\bpr

From \eqref{ED2}, we have
\[
 \int_{\Omega}  \frac{1}{|\omega_{\epsilon}|}\chi_{\omega_{\epsilon}}(x) \frac{\partial w_{\epsilon}^i}{\partial x_j}(x) \phi(x) d x \rightarrow \int_{\Omega} \phi(x)M^{\theta}_{ij} d \mu^{\theta}\quad  \mbox{ for all } \phi \in C^0(\overline{\Omega}).
\]
Multiplying both sides by $(\gamma_1-\gamma_0)$, we get
\[
 \int_{\Omega} (\gamma_1-\gamma_0) \frac{1}{|\omega_{\epsilon}|}\chi_{\omega_{\epsilon}}(x) \frac{\partial w_{\epsilon}^i}{\partial x_j}(x) \phi(x) d x \rightarrow \int_{\Omega} (\gamma_1-\gamma_0)\phi(x)M^{\theta}_{ij} d\mu^{\theta} \quad\mbox{ for all } \phi \in C^0(\overline{\Omega}).
 \]
Since  $\gamma_{\epsilon}(x) = \chi_{\omega_{\epsilon}}(x)\gamma_1 + (1-\chi_{\omega_{\epsilon}}(x))\gamma_0$, we have,
\[
\int_{\Omega} (\gamma_1-\gamma_0)\phi(x)M^{\theta}_{ij} d \mu^{\theta} = \lim_{\epsilon\to 0} \frac{1}{|\omega_{\epsilon}|}\int_{\Omega} (\gamma_{\epsilon}(x) - \gamma_0)\frac{\partial w_{\epsilon}^i}{\partial x_j} \phi(x) d x.
\]
From \eqref{Corrector Theory},  
\[
\int_{\Omega} (\gamma_1-\gamma_0)\phi(x)M^{\theta}_{ij} d \mu^{\theta} =\frac{1}{\int\theta(x)d x}\int\limits_{\Omega} (\gamma^{*}_{ij}(x)-\gamma_0\delta_{ij})\phi(x) d x.
\]
Now since $ d \mu^{\theta} = \frac{\theta(x)}{\int_{\Omega} \theta(x)d x}d x \ $ almost everywhere in $\Omega$
we obtain \eqref{PT and HT relation}. This completes the proof.
\hfill\epr
\begin{remark}[Localization principle for the polarization tensor of the non zero volume fraction]
Let $\gamma_\epsilon\in M(\gamma_1,\gamma_0;\Omega)$ and $\widetilde{\gamma}_\epsilon\in M(\gamma_1,\gamma_0;\Omega)$ 
be two sequences, which H-converge to $\gamma^{*}(x)$ and $\widetilde{\gamma}^{*}(x)$, 
respectively. Let $U$ be an open subset compactly embedded in $\Omega$ and if
$\gamma_\epsilon(x)=\widetilde{\gamma}_\epsilon(x)$ in $U$. It is known that 
$\gamma^{*}(x)=\widetilde{\gamma}^{*}(x)$ in $U$ with the same volume fraction $\theta$ in $U$.
Thus it follows from \eqref{PT and HT relation} $M^{\theta}(x)=\widetilde{M}^{\theta}(x)$ in $U$. 
\end{remark}
\begin{remark}
We don't have a similar property for polarization tensor corresponding to near
zero volume fraction. 
Let $\gamma_\epsilon\in M(\gamma_1,\gamma_0;\Omega)$ and $\widetilde{\gamma}_\epsilon\in M(\gamma_1,\gamma_0;\Omega)$ 
be two sequences with zero- volume fractions i.e. $|\omega_\epsilon|$ and $|\widetilde{\omega}_\epsilon|$ goes to zero
as $\epsilon$ tends to zero. Let $U$ be an open subset compactly embedded in $\Omega$ and 
$\gamma_\epsilon(x)=\widetilde{\gamma}_\epsilon(x)$ in $U$. Then if in addition  
$\frac{|\omega_\epsilon|}{|\widetilde{\omega}_\epsilon|} \rightarrow 1$ as $\epsilon \rightarrow 0$,
 it follows from \eqref{m0-d} and \eqref{mu0-d} that
$M^{0}(x)=\widetilde{M}^{0}(x)$ and $\mu^{0}(x)=\widetilde{\mu}^{0}(x)$ in $U$.
\end{remark}
\begin{remark} Using Proposition \ref{PT and HT Proposition} and \eqref{ED5}, we can rewrite the asymptotic formula in Theorem \ref{Polarization tensor thm-a} as follows: 
$$
(\gamma_{\epsilon}\frac{\partial u_{\epsilon}}{\partial \nu} - \gamma^{*}\frac{\partial u}{\partial \nu})(y) 
=\ (|\omega_{\epsilon}|-\delta)\int_{\Omega} (\gamma_1-\gamma_0) M^{\theta}_{ij}(x)\frac{\partial u}{\partial x_i}(x)\frac{\partial D}{\partial x_j}(x,y)d \mu^{\theta}(x) + o(1).
$$
\end{remark}
\begin{corollary}
The polarization tensor $M^{\theta}$ is symmetric.
\end{corollary}
\bpr
This follows from \eqref{PT and HT relation} and from the fact that $\gamma^{*}$ is symmetric. 
\hfill\epr
\subsection{Bounds on $M^{\theta}$}
We now derive bounds for the polarization tensor based on relation we obtained between the this tensor and the homogenization tensor (Proposition \ref{PT and HT Proposition}).
\begin{proposition}\label{bound-mt} We have the following bounds for the polarization tensor $M^\theta$ for $\mu^\theta$ almost everywhere $x$
\[ \mathrm{min}\left\{ 1, \frac{\gamma_0}{\theta \gamma_0 +(1-\theta)\gamma_1} \right\}|\xi|^2 \leq M^{\theta}_{ij}(x)\xi_i\xi_j \leq \mathrm{max}\left\{ 1, \frac{\gamma_0}{\theta \gamma_0 +(1-\theta)\gamma_1} \right\}|\xi|^2 , \hspace{4pt} \xi \in \mathbb{R}^N. 
\]
\end{proposition}
\bpr
The proof is a straightforward application of Proposition \ref{PT and HT Proposition} the arithmetic mean $(\gamma_aI$), harmonic mean $(\gamma_hI$) bound for $\gamma^{*}$ \cite[\S 2.1.2]{A}: 
\[
\underbrace{\lb\frac{\theta(x)}{\gamma_1}+\frac{1-\theta(x)}{\gamma_0}\rb^{-1}}_{\gamma_h}I\  \leq\  \gamma^{*}(x)\  \leq  \ \underbrace{(\theta(x) \gamma_1+(1-\theta(x))\gamma_0)}_{\gamma_a}I.
\]
\hfill\epr
The above result generalizes to the case of non-zero volume fraction of the conductivities, the bounds obtained at the conclusion of \cite[Lemma 3]{CV1}. 

\begin{proposition}\label{Trace bounds - general} 
We have the following pointwise trace bounds for the polarization tensor $M^{\theta}$ for $\mu^\theta$ almost everywhere $x$
\begin{align}
\mbox{Upper bound}:\quad &\mathrm{trace}\lb I-\theta(x) M^{\theta}(x)\rb^{-1}\leq \frac{N}{1-\theta(x)} + \frac{\theta(x)}{1-\theta(x)}\Big(\frac{\gamma_0}{\gamma_1} - 1\Big),\ \label{ub} \\
\mbox{Lower bound}:\quad &\mathrm{trace}\lb \theta(x) M^{\theta}(x)\rb ^{-1} \leq \frac{N}{\theta(x)} - \frac{1-\theta(x)}{\theta(x)}\Big( 1 - \frac{\gamma_1}{\gamma_0}\Big).\ \label{lb}
\end{align}
\end{proposition}
\bpr The proof follows by slight modifications of the proof of \cite [Prop. 3.1]{AM}, \cite{T} and the relation between the polarization tensor and the homogenization tensor ( cf. Proposition \ref{PT and HT Proposition} ).
Recall that we have let $\gamma_\epsilon(x) = \chi_{\epsilon}(x)\gamma_1 +(1 - \chi_\epsilon(x))\gamma_0$.
Following the proof as in \cite[Prop. 3.1]{AM}, we have the following 
pointwise lower bound for $\gamma^{*}(x)$
\begin{equation}\label{ED4}
\lb\gamma^{*}(x)-\gamma_1I\rb^{-1}\leq \frac{I}{\gamma_{a}(x)-\gamma_1}+\frac{\theta(x)(1-\theta(x))(\gamma_1-\gamma_0)^{2}}{\gamma_1(\gamma_{a}(x)-\gamma_1)^{2}}\widetilde{M}(x),
\end{equation}
where $\widetilde{M}$ \cite[Eq. 12]{AM} is given by 
$$\widetilde{M}(x) = \frac{1}{\theta(x)(1-\theta(x))}\int_{\mathbb{S}^{N-1}}\xi\otimes\xi\ d\nu(\xi).$$
Here $\theta(x)(1-\theta(x))dx\otimes d\nu(\xi)$ is the $H$-measure ( \cite[\S 28]{T} ) of the sequence $(\chi_{\o_{\epsilon}}-\theta)(x)$. Note that $\widetilde{M}$ has unit trace.

Now from Proposition \ref{PT and HT Proposition}, we have 
\[
\gamma^{*}(x)-\gamma_0I=(\gamma_1-\gamma_0)\theta(x)M^{\theta}(x).
\]
Now substituting this into \eqref{ED4}, we have
\[
\lb I-\theta(x)M^{\theta}(x)\rb^{-1}\ \leq\ \frac{\gamma_0-\gamma_1}{\gamma_{a}(x)-\gamma_1} I + \frac{\theta(x)(1-\theta(x))(\gamma_0-\gamma_1)^{3}}{\gamma_1(\gamma_{a}(x)-\gamma_1)^{2}}\widetilde{M}(x)
\]
Since $\gamma_{a}(x)-\gamma_1=(1-\theta(x))(\gamma_0-\gamma_1)$ and then by taking trace, 
we get the inequality \eqref{ub}.

Similarly, by using upper bound for $\gamma^{*}$ \cite[Prop. 3.1]{AM}, we deduce \eqref{lb} which gives the lower
bound for $M^{\theta}$. 
\hfill\epr
\begin{remark}
Proposition \ref{PT and HT Proposition} combined with \cite[Theorem 2.2.13] {A} gives us that the pointwise trace bounds for the polarization tensor $M^{\theta}$ obtained in the previous proposition are optimal.
\end{remark}

\section{Proof of Theorem \ref{approximation} and Theorem \ref{Optimal trace bounds theorem}}
\setcounter{equation}{0}

\bpr[Proof of Theorem \ref{approximation}]
\textbf{Step 1 : }
Given $d\mu^{0}$, $M^{0}$ with $M^0 \in L^2(\Omega,d \mu^{0})$, there exists a sequence of 
microstructures $\omega_{\epsilon} \subset K \subset\Omega$ (for some compact set $K$, as $dist(\omega_{\epsilon},\partial\Omega)\geq d_0 >0$) such that as $\epsilon \rightarrow 0, \hspace{3pt}|\omega_\epsilon| \rightarrow 0$ and
\begin{equation}\begin{aligned}\label{qq7}
(i)\ \ \frac{1}{|\omega_{\epsilon}|}\chi_{\omega_{\epsilon}}(x)& \stackrel{*}\rightharpoonup d \mu^{0} \mbox{ in } (C^0(\overline{\Omega}))^{*}.\\
(ii)\ \ \frac{1}{|\omega_{\epsilon}|}\chi_{\omega_{\epsilon}}(x) \frac{\partial w_{\epsilon}^i}{\partial x_j}dx &\stackrel{*}\rightharpoonup d M^{0}_{ij}= M^{0}_{ij}d \mu^{0}  \textrm{ in } (C^0(\overline{\Omega}))^{*}.
\end{aligned}\end{equation}
where $w_{\epsilon}^{i}$ are correctors defined in \eqref{ED6}.\\
\\
We consider any point $x_0 \in Support\ of\ (d\mu^{0})$ and an open
cube $Q_{x_0, h} = x_0 + (-\frac{h}{2},\frac{h}{2})^N$ centered at 
the point $x_0$ which is included in $K$ for sufficiently small $h>0$. 
Note we can write $Q_{x_0, h} = x_0 + hY$, 
where the cell $Y = (-\frac{1}{2},\frac{1}{2})^N$ with $|Y|=1$. 
As $x_0 \in Support\ of\ (d\mu^{0})$, 
\begin{center}$\ \mu^{0}(Q_{x_0,h}) >0$,\ $\forall h>0.$\end{center}
Making the change of variables
$$x\in Q_{x_0,h} \mapsto y\in Y\mbox{ as }x= x_0 +hy$$
we localize  $d\mu^0$ and $d M^0_{ij}$ over $Q_{x_0,h}$.\\
\\
Now following the idea that locally any inhomogeneous medium can be approximated by the periodic medium 
( cf. \cite[Theorem 1.3.23]{A} ), let us consider the periodic homogenization with inhomogeneities
$(\omega_\epsilon\cap Q_{x_0,h})$ contained in $Q_{x_0,h}$ and extending it periodically in whole $\mathbb{R}^N$ with the small period $\eta>0$. 
The corresponding coefficient is denoted by $\gamma_{\epsilon}(x_0 +h\frac{x}{\eta})$. 
Let $\gamma^{*}_{x_0, \epsilon, h}$ denote the homogenized tensor thus obtained. 
The volume fraction of inhomogeneities is evidently 
\begin{equation}\label{qq9}\theta_{x_0,\epsilon,h} = \int_Y \chi_{\omega_\epsilon \cap Q_{x_0,h}}(x_0 +hy)\ dy=\ \frac{|\omega_\epsilon\cap Q_{x_0,h}|}{|Q_{x_0,h}|}.\end{equation}
The volume fraction $\theta_{x_0,\epsilon,h} > 0$ for $\epsilon$ being small enough,
which follows from \eqref{qq3} below.
We recall the integral representation of $\gamma^{*}_{x_0,\epsilon,h}$ 
\begin{equation}\label{gst} (\gamma^{*}_{x_0, \epsilon, h})_{ij} = \int_{Y} \gamma_{\epsilon}(x_0 +hy)\nabla\widetilde{w}^{i}_{x_0,\epsilon,h}(y)\cdot e_j\ dy\end{equation}
where, $(\widetilde{w}^{i}_{x_0,\epsilon,h}(y))_{1\leq i\leq N}$ is the family of unique solutions in $H^{1}(Y)/{\mathbb{R}}$ of the cell problems
\begin{equation}\label{w-t-c}\nabla\cdot(\gamma_{\epsilon}(x_0 + hy)\nabla\widetilde{w}^{i}_{x_0,\epsilon,h}(y)) = 0 \quad\mbox{in }Y,\quad y \mapsto (\widetilde{w}^{i}_{x_0,\epsilon,h}(y)-y_i) \quad\mbox{ is } Y \mbox{ periodic.}\end{equation}
In the sequel, we consider the polarization tensor denoted as $M^{\theta_{x_0,\epsilon,h}}(x_0)$
which corresponds to the above periodic microstructure. As observed in Remark \ref{per-cons}, $M^{\theta_{x_0,\epsilon,h}}(x_0)$
is a constant matrix. Using the relation \eqref{PT and HT relation} between the homogenized tensor and the polarization tensor
with non-zero volume fraction $\theta_{x_0,\epsilon,h},$ the following integral representation is easily obtained
from \eqref{gst} :
\begin{align}
M_{ij}^{\theta_{x_0, \epsilon, h}}(x_0)&=\ \int_{Y}\frac{1}{\theta_{x_0,\epsilon, h}}\chi_{\omega_{\epsilon}\cap Q_{x_0,h} }(x_0+ hy) \nabla_y\widetilde{w}^{i}_{x_0,\epsilon,h}(y)\cdot e_{j}\ dy\notag\\
                                       &=\ \frac{1}{|\omega_\epsilon\cap Q_{x_0,h}|}\int_{Q_{x_0,h}}\chi_{\omega_{\epsilon}\cap Q_{x_0,h} }(x)\ h\nabla_x \widetilde{w}^i_{x_0,\epsilon,h}(\frac{x-x_0}{h})\cdot e_{j}\ dx.\label{mst}
\end{align}
Our next task is to replace $\widetilde{w}^{i}_{x_0,\epsilon,h}$ by $w^i_\epsilon$ in \eqref{mst} as well as to analyze its limiting behavior as $\epsilon \rightarrow 0$ and $h \rightarrow 0$ in that order.   
In order to do that, we invoke our next step as follows. \\
\\
\textbf{Step 2 : }
We begin with defining limiting quantities representing volume fraction of inhomogeneities in $Q_{x_0,h}$ $(h>0):$ 
\begin{align*}
 \underset{\epsilon \rightarrow 0}{lim\hspace{1pt}sup}\ \frac{|\omega_{\epsilon}\cap Q_{x_0,h}|}{|\omega_{\epsilon}|}\ =\ V^{+}_h(x_0) \ \ \mbox{ and } \ 
\underset{\epsilon \rightarrow 0}{lim\hspace{1pt}inf}\ \frac{|\omega_{\epsilon}\cap Q_{x_0,h}|}{|\omega_{\epsilon}|}\ =\ V^{-}_h(x_0) \ \  
 \end{align*} 
Clearly
$V^{+}_h(x_0)$ and $V^{-}_h(x_0)$ are monotonically increasing function with respect to $h$, i.e.
\begin{center} $ V^{+}_h(x_0) \ \ ( \mbox{ respectively, } V^{-}_h(x_0)\ )\ \geq\ V^{+}_{\widetilde{h}}(x_0)\ \ ( \mbox{ respectively, } V^{-}_{\widetilde{h}}(x_0)\ ),\ \mbox{ whenever }h \geq \widetilde{h}$.\end{center}
\textbf{ Claim : } For $h>0$ fixed,
\begin{equation}\label{qq3} \mu^{0}(\overline{Q}_{x_0,h})\ =\ V_h^{+}(x_0)\ =\ V_h^{-}(x_0). \end{equation}
\begin{remark}\label{muR}
 It follows that 
$\underset{\epsilon \rightarrow 0}{lim}\ \frac{|\omega_{\epsilon}\cap Q_{x_0,h}|}{|\omega_{\epsilon}|}$ exists and equal to $\mu^{0}(\overline{Q}_{x_0,h}) > 0$.  
\end{remark}
\textit{Proof of the above claim : }
We will first show 
\begin{center}$V^{+}_h(x_0) \leq \mu^{0}(\overline{Q}_{x_0,h})$.\end{center}
Let us consider an open cube $Q_{x_0,h-\delta}$ for $\delta >0$ small enough.
Note $Q_{x_0,h-\delta}$ is relatively compact in $Q_{x_0,h}$.
We consider a sequence of test functions $\phi^\delta_h\in C^{0}(\overline{\Omega})$ for $h$ fixed such that
\begin{equation}\label{qq6} Support\ of\ \phi^\delta_h \subset Q_{x_0,h},\ \phi^\delta_h\equiv 1 \ \mbox{ in }\ Q_{x_0,h-\delta}, \ \mbox{ and }\ 0\leq \phi^\delta_h(x) \leq 1,\ \ x\in  Q_{x_0,h}\smallsetminus \overline{Q}_{x_0,h-\delta}.\end{equation}   
As we see for $h>0$ fixed , as $\delta \downarrow 0$,
\begin{center}$ \phi^\delta_h(x)\ \rightarrow\ 1 \ \mbox{  everywhere $x$ in } \ Q_{x_0,h}.$ \end{center}
Thus using the Lebesgue dominated convergence theorem and the fact $\mu^{0}$ is a bounded measure
$$ \mu^{0}(Q_{x_0,h}) = \int_{Q_{x_0,h}} d\mu^{0} =\ \underset{\delta\downarrow 0}{lim}\ \int_{Q_{x_0,h}} \phi^\delta_h d\mu^{0}. $$
So,
\begin{align}
\mu^{0}(Q_{x_0,h})=\ \underset{\delta\downarrow 0}{lim}\ \int_{Q_{x_0,h}} \phi^\delta_h d\mu^{0}=\ \underset{\delta\downarrow 0}{lim}\ \int_{\Omega} \phi^\delta_h d\mu^{0}
&=\  \underset{\delta\downarrow 0}{lim}\lb \underset{\epsilon\rightarrow 0}{lim}\ \int_{\Omega}\frac{1}{|\omega_{\epsilon}|}\chi_{\omega_{\epsilon}}\phi^\delta_h dx\rb,\notag\\
&\geq\ \underset{\delta\downarrow 0}{lim\hspace{1pt}sup}\lb\underset{\epsilon\rightarrow 0}{lim\hspace{1pt}sup}\ \int_{Q_{x_0,h-\delta}}\frac{1}{|\omega_{\epsilon}|}\chi_{\omega_{\epsilon}}\phi^\delta_h dx\rb,\notag\\
&=\  \underset{\delta\downarrow 0}{lim\hspace{1pt}sup}\lb\underset{\epsilon\rightarrow 0}{lim\hspace{1pt}sup}\ \int_{Q_{x_0,h-\delta}}\frac{1}{|\omega_{\epsilon}|}\chi_{\omega_{\epsilon}} dx\rb, \notag\\
&=\  \underset{\delta\downarrow 0}{lim\hspace{1pt}sup}\lb\underset{\epsilon\rightarrow 0}{lim\hspace{1pt}sup}\ \frac{|\omega_{\epsilon}\cap Q_{x_0,h-\delta}|}{|\omega_{\epsilon}|}\rb,\notag\\
&=\  \underset{\delta\downarrow 0}{lim\hspace{1pt}sup}\ V^{+}_{h-\delta}(x_0).\label{qq1}
\end{align}
Now as $\mu^{0}$ is a Radon measure, for any given $\eta >0$ there exist a $\delta_\eta>0$ such that
$$ \mu^{0}(\overline{Q}_{x_0,h}) +\eta\ \geq \ \mu^{0}(Q_{x_0,h+\delta_\eta}) \quad\mbox{ ( outer regularity )}.$$
And using the above result \eqref{qq1} with $(h+\delta_\eta)$ in the place of $h$, we get 
$$ \mu^{0}(\overline{Q}_{x_0,h}) +\eta\ \geq \ \underset{\delta\downarrow 0}{lim\hspace{1pt}sup}\ V^{+}_{h+\delta_\eta-\delta}(x_0).$$
Now as $0 < \delta < \delta_\eta $, by using the monotonicity we get
$ V^{+}_{h+\delta_\eta -\delta}(x_0) \geq V^{+}_h(x_0)$ for $h$ being fixed. Consequently,
$$\mu^{0}(\overline{Q}_{x_0,h}) +\eta\ \geq \ V^{+}_{h}(x_0).$$
As $\eta > 0$ is chosen arbitrarily therefore  
$$\mu^{0}(\overline{Q}_{x_0,h})\ \geq \ V^{+}_{h}(x_0).$$
\\
Similarly, we will show 
$$ \mu^{0}(\overline{Q}_{x_0,h})\ \leq \ V^{-}_{h}(x_0).$$
Since the arguments are slightly different, we go through them.
Here we will consider an open cube $Q_{x_0,h+\delta}$ for $\delta >0$ small enough. Note that
$Q_{x_0,h}$ is relatively compact in $Q_{x_0,h+\delta}$ for every $\delta >0$. 
We consider a sequence of test functions $\psi^\delta_h\in C^{0}(\overline{\Omega})$ for $h$ fixed, such that
\begin{equation}\label{pshd} Support\ of\ \psi^\delta_h\subset Q_{x_0,h+\delta},\  \psi^\delta_h\equiv 1 \ \mbox{ in }\ \overline{Q}_{x_0,h}, \ \mbox{ and }\ 0\leq \psi^\delta_h(x) \leq 1,\ \ x\in  Q_{x_0,h+\delta}\smallsetminus \overline{Q}_{x_0,h}.\end{equation}   
So,
$$\lb \int_{\Omega}\psi^\delta_h d\mu^{0} - \int_{\overline{Q}_{x_0,h}}d\mu^{0}\rb = \int_{Q_{x_0,h+\delta}\smallsetminus \overline{Q}_{x_0,h}}\psi^\delta_h d\mu^{0} \leq\ \mu^{0}(Q_{x_0,h+\delta}-\overline{Q}_{x_0,h}) \rightarrow 0 \ \ \mbox{ as }\delta\rightarrow 0.$$ 
Thus 
\begin{align}
\mu^{0}(\overline{Q}_{x_0,h}) = \int_{\overline{Q}_{x_0,h}} d\mu^{0} =\ \underset{\delta\downarrow 0}{lim}\ \int_{\Omega} \psi^\delta_h d\mu^{0}&=\  \underset{\delta\downarrow 0}{lim}\lb \underset{\epsilon\rightarrow 0}{lim}\ \int_{\Omega}\frac{1}{|\omega_{\epsilon}|}\chi_{\omega_{\epsilon}}\psi^\delta_h dx\rb\notag\\
&\leq\ \underset{\delta\downarrow 0}{lim\hspace{1pt}inf}\lb\underset{\epsilon\rightarrow 0}{lim\hspace{1pt}inf}\ \int_{Q_{x_0,h+\delta}}\frac{1}{|\omega_{\epsilon}|}\chi_{\omega_{\epsilon}} dx\rb\notag\\
&=\  \underset{\delta\downarrow 0}{lim\hspace{1pt}inf}\lb\underset{\epsilon\rightarrow 0}{lim\hspace{1pt}inf}\ \frac{|\omega_{\epsilon}\cap Q_{x_0,h+\delta}|}{|\omega_{\epsilon}|}\rb\notag\\
&=\  \underset{\delta\downarrow 0}{lim\hspace{1pt}inf}\ V^{-}_{h+\delta}(x_0)\label{qq2}.
\end{align}
Now as $\mu^{0}$ is a Radon measure, for any given $\eta >0$ there exists a $\delta_\eta>0$ such that
$$ \mu^{0}(\overline{Q}_{x_0,h}) -\eta\ \leq \ \mu^{0}(\overline{Q}_{x_0,h-\delta_\eta}) \quad\mbox{ ( inner regularity )}.$$
And using the above result \eqref{qq2} with $(h-\delta_\eta)$ in the place of $h$, we get 
$$ \mu^{0}(\overline{Q}_{x_0,h}) -\eta\ \leq \ \underset{\delta\downarrow 0}{lim\hspace{1pt}inf}\ V^{-}_{h-\delta_\eta+\delta}(x_0).$$
Now as $0 < \delta < \delta_\eta $, by using the monotonicity we get
$ V^{-}_{h-\delta_\eta +\delta}(x_0) \leq V^{-}_h(x_0)$ for $h$ being fixed. It follows
$$\mu^{0}(\overline{Q}_{x_0,h}) -\eta\ \leq \ V^{-}_{h}(x_0).$$
As $\eta > 0$ is chosen arbitrarily therefore  
$$\mu^{0}(\overline{Q}_{x_0,h})\ \leq \ V^{-}_{h}(x_0).$$
Hence we have established our claim \eqref{qq3}.\\
\\
\textbf{ Step 3 : }
As the next step, in a very similar way as we just did in Step 2, we will show, for $h>0$ fixed, that 
\begin{equation}\label{qq4} \int_{\overline{Q}_{x_0,h}} d M^{0}_{ij}(x) \ =\ V^{+}_{h}(x_0)\cdot \underset{\epsilon\rightarrow 0}{lim\hspace{1pt}sup}\ M^{\theta_{x_0,\epsilon,h}}_{ij}(x_0) =\  V^{-}_{h}(x_0)\cdot \underset{\epsilon\rightarrow 0}{lim\hspace{1pt}inf}\ M^{\theta_{x_0,\epsilon,h}}_{ij}(x_0).\end{equation}
Let us first show that, 
$$\int_{\overline{Q}_{x_0,h}} d M^{0}_{ij}(x) \ \geq\ V^{+}_{h}(x_0)\cdot \underset{\epsilon\rightarrow 0}{lim\hspace{1pt}sup}\ M^{\theta_{x_0,\epsilon,h}}_{ij}(x_0) $$
Using the Lebesgue dominated convergence theorem together with the fact $d M^{0}_{ij}$( cf. \eqref{qq7} ) is a bounded measure we write first
\begin{align}
\int_{Q_{x_0,h}} d M^{0}_{ij}(x) =\  \underset{\delta\downarrow 0}{lim}\ \int_{Q_{x_0,h}} \phi^\delta_h d M^{0}_{ij}&=\ \underset{\delta\downarrow 0}{lim}\ \int_{\Omega} \phi^\delta_h d M^{0}_{ij}\notag\\
&=\ \underset{\delta\downarrow 0}{lim}\lb \underset{\epsilon\rightarrow 0}{lim}\ \int_{\Omega}\frac{1}{|\omega_{\epsilon}|}\chi_{\omega_{\epsilon}}(x)\nabla_x w^i_{\epsilon}(x)\cdot e_{j}\ \phi^\delta_h(x)\ dx\rb\label{mwd}
\end{align}
where, $\phi^\delta_h$ is defined as in \eqref{qq6}.\\
\\
Our next tusk is to replace $\omega_\epsilon$ by $\widetilde{\omega}_{x_0,\epsilon,h}$ in the above representation.\\
\\
\textbf{Claim : }  
\begin{align}
&\frac{1}{|\omega_\epsilon\cap Q_{x_0,h}|}\int_{Q_{x_0,h}}\chi_{\omega_{\epsilon}\cap Q_{x_0,h} }(x)\nabla_x {w }^{i}_{\epsilon}(x)\cdot e_{j}\ \phi^\delta_h(x)\ dx \notag\\
&\quad\quad\quad\quad\quad=\ \frac{1}{|\omega_\epsilon\cap Q_{x_0,h}|} \int_{Q_{x_0,h}}\chi_{\omega_{\epsilon}\cap Q_{x_0,h} }(x)\ h\nabla_x\widetilde{w}^i_{x_0,\epsilon,h}(\frac{x-x_0}{h})\cdot e_{j}\ \phi^\delta_h(x)\ dx+ \mathcal{E}_{h,\phi^\delta_h}(\epsilon)\notag\\
&\mbox{where, }\ \mathcal{E}_{h,\phi^\delta_h}(\epsilon)\rightarrow 0\  \mbox{as }\epsilon \rightarrow 0\  \mbox{ for every fixed }h \mbox{ and }\phi^\delta_h.\label{dhe}
\end{align}
As we see, from \eqref{w-epsilon} and \eqref{w-t-c} it follows that
\begin{align*}
&\gamma_1||\phi^\delta_h\lb\nabla_x w^i_\epsilon(x)-h\nabla_x\widetilde{w}^i_{x_0,\epsilon,h}(\frac{x-x_0}{h})\rb||^2_{L^2(Q_{x_0,h})}\\
&\leq \int_{Q_{x_0,h}}(\phi^\delta_h)^2\gamma_\epsilon(x)\lb\nabla_x w^i_\epsilon(x)-h\nabla_x\widetilde{w}^i_{x_0,\epsilon,h}(\frac{x-x_0}{h})\rb\cdot\lb\nabla_x w^i_\epsilon(x)-h\nabla_x\widetilde{w}^i_{x_0,\epsilon,h}(\frac{x-x_0}{h})\rb dx\\
&= -\int_{Q_{x_0,h}}\lb \nabla_x\cdot\gamma_\epsilon(x)\nabla_x w^i_\epsilon(x)-\nabla_x\cdot\gamma_\epsilon(x)h\nabla_x\widetilde{w}^i_{x_0,\epsilon,h}(\frac{x-x_0}{h})\rb(\phi^\delta_h)^2\lb w^i_\epsilon(x)-h\ \widetilde{w}^i_{x_0,\epsilon,h}(\frac{x-x_0}{h})\rb dx \\
&\quad-\int_{Q_{x_0,h}}\lb w^i_\epsilon(x)-h\ \widetilde{w}^i_{x_0,\epsilon,h}(\frac{x-x_0}{h})\rb\gamma_\epsilon(x)\lb\nabla_x w^i_\epsilon(x)-\nabla_x\widetilde{w}^i_{x_0,\epsilon,h}(\frac{x-x_0}{h})\rb\cdot\nabla_x(\phi^\delta_h)^2dx\\
&= -\int_{Q_{x_0,h}}\gamma_\epsilon(x)\lb w^i_\epsilon(x)-h\ \widetilde{w}^i_{x_0,\epsilon,h}(\frac{x-x_0}{h})\rb2\nabla_x\phi^\delta_h\cdot \phi^\delta_h\lb \nabla_x w^i_\epsilon(x)-h\nabla_x\widetilde{w}^i_{x_0,\epsilon,h}(\frac{x-x_0}{h})\rb dx.
\end{align*}
Thus,
\begin{equation*}\begin{aligned}
C\ ||\phi^\delta_h\lb\nabla_x w^i_\epsilon(x)-h\nabla_x\widetilde{w}^i_{x_0,\epsilon,h}(\frac{x-x_0}{h})\rb&||_{L^2(Q_{x_0,h})}\\
\leq \ ||w^i_\epsilon(x)&-x_i||_{L^2(Q_{x_0,h})}+||h\ \widetilde{w}^i_{x_0,\epsilon,h}(\frac{x-x_0}{h})-x_i||_{L^2(Q_{x_0,h})}.
\end{aligned}\end{equation*}
Now, as it is shown in \cite[Lemma 1.]{CV1}, we invoke 
\begin{equation}\label{cvst}
\begin{aligned}
&||w^i_\epsilon(x)-x_i||_{L^2(\Omega)} \leq\ o(|\omega_\epsilon|^{\frac{1}{2}}) =\ o(|\omega_\epsilon\cap Q_{x_0,h}|^{\frac{1}{2}})\ \ \mbox{ (by using the Remark }\ref{muR}) \\
\mbox{and}\ \ &||h\ \widetilde{w}^i_{x_0,\epsilon,h}(\frac{x-x_0}{h})-x_i||_{L^2(Q_{x_0,h})}\leq\ o(|\omega_\epsilon\cap Q_{x_0,h}|^{\frac{1}{2}}).
\end{aligned}
\end{equation}
Hence we get
\begin{align*}
|\mathcal{E}_{h,\phi^\delta_h}(\epsilon)|&\leq\ \frac{1}{|\omega_\epsilon\cap Q_{x_0,h}|}\ ||\chi_{\omega_\epsilon\cap Q_{x_0,h}}(x)\ \phi^\delta_h\lb\nabla_x w^i_\epsilon(x)-h\nabla_x\widetilde{w}^i_{x_0,\epsilon,h}(\frac{x-x_0}{h})\rb||_{L^1(Q_{x_0,h})}\\
&\leq \frac{1}{|\omega_\epsilon\cap Q_{x_0,h}|^{\frac{1}{2}}}\ ||\phi^\delta_h\lb \nabla_x w^i_\epsilon(x)-h\nabla_x\widetilde{w}^i_{x_0,\epsilon,h}(\frac{x-x_0}{h})\rb||_{L^2(Q_{x_0,h})}\\
& =\ o(1) \rightarrow 0 \quad\mbox{as }\epsilon\rightarrow 0 \ \mbox{ ( by using \eqref{cvst} for being $h$,$\phi^\delta_h$ fixed ). }
\end{align*}
\\
Thus from \eqref{mwd} and \eqref{dhe} we get
\begin{align}
&\int_{Q_{x_0,h}} d M^{0}_{ij}(x) =\ \underset{\delta\downarrow 0}{lim}\lb \underset{\epsilon\rightarrow 0}{lim}\ \int_{\Omega}\frac{1}{|\omega_{\epsilon}|}\chi_{\omega_{\epsilon}}(x)\nabla_x {w }^{i}_{\epsilon}(x)\cdot e_{j}\ \phi^\delta_h(x)\ dx\rb\notag \\
&\quad =\ \underset{\delta\downarrow 0}{lim}\lb \underset{\epsilon \rightarrow 0}{lim}\ \frac{|\omega_{\epsilon}\cap Q_{x_0,h}|}{|\omega_{\epsilon}|}\ \frac{1}{|\omega_\epsilon\cap Q_{x_0,h}|}\int_{Q_{x_0,h}}\chi_{\omega_{\epsilon} }(x)\ h\nabla_x \widetilde{w}^{i}_{\epsilon}(\frac{x-x_0}{h})\cdot e_{j}\ \phi^\delta_h(x)\ dx\rb. \label{qq5} 
\end{align}
We show below the contribution coming from the annular region $Q_{x_0,h}\smallsetminus Q_{x_0,h-\delta}$ is negligible. We may also replace $|\omega_\epsilon\cap Q_{x_0,h}|$ by 
$|\omega_\epsilon\cap Q_{x_0,h-\delta}|$ as shown below.
\begin{align}
&\frac{1}{|\omega_\epsilon\cap Q_{x_0,h}|} \int_{Q_{x_0,h}\smallsetminus Q_{x_0,h-\delta}}\chi_{\omega_\epsilon}(x)\ h\nabla_x \widetilde{w}^{i}_{\epsilon}(\frac{x-x_0}{h})\cdot e_j\ \phi^{\delta}_h(x)\ dx \notag\\
&\frac{1}{|\omega_\epsilon\cap Q_{x_0,h}|} \int_{Q_{x_0,h}\smallsetminus Q_{x_0,h-\delta}}\chi_{\omega_\epsilon}(x)\lb h\nabla_x \widetilde{w}^{i}_{\epsilon}(\frac{x-x_0}{h})-e_i\hspace{1.2pt}+e_i\rb\cdot e_j\ \phi^{\delta}_h(x)\ dx \notag\\
&\leq\ C\ \frac{1}{|\omega_\epsilon\cap Q_{x_0,h}|}\lb |\omega_\epsilon\cap Q_{x_0,h}\smallsetminus Q_{x_0,h-\delta}|^{\frac{1}{2}}\rb.\lb|\omega_\epsilon\cap Q_{x_0,h}|^{\frac{1}{2}} + |\omega_\epsilon \cap Q_{x_0,h}\smallsetminus Q_{x_0,h-\delta}|^{\frac{1}{2}}\rb\notag\\
&\qquad\qquad\qquad\qquad\qquad\qquad\qquad\qquad\qquad\qquad\qquad\qquad\mbox{( using \eqref{ED1}, cf. Remark \eqref{qq8} )}\notag\\
&=\ C \lb (1- \frac{|\omega_\epsilon\cap Q_{x_0,h-\delta}|}{|\omega_\epsilon\cap Q_{x_0,h}|})^{\frac{1}{2}} + ( 1- \frac{|\omega_\epsilon\cap Q_{x_0,h-\delta}|}{|\omega_\epsilon\cap Q_{x_0,h}|})\rb. \label{hder}
\end{align}
Now as we see that, 
\begin{align}\label{depl}
\underset{\delta\downarrow 0}{lim}\lb\underset{\epsilon\rightarrow 0}{lim}\ \frac{|\omega_\epsilon\cap Q_{x_0,h-\delta}|}{|\omega_\epsilon\cap Q_{x_0,h}|}\rb 
=\ \underset{\delta\downarrow 0}{lim}\lb\underset{\epsilon\rightarrow 0}{lim}\ \frac{|\omega_{\epsilon}\cap Q_{x_0,h}|/|\omega_\epsilon|}{|\omega_{\epsilon}\cap Q_{x_0,h-\delta}|/|\omega_\epsilon|}\rb =\ \underset{\delta\downarrow 0}{lim}\lb\frac{\mu^{0}(\overline{Q}_{x_0,h})}{\mu^{0}(\overline{Q}_{x_0,h-\delta})}\rb=\ 1.&\\
(\because \mu^0 \mbox{ is a Radon measure })&\notag
\end{align}
So, \eqref{hder} goes to zero as $\epsilon \rightarrow 0$ and $\delta \downarrow 0$ in that order, for every fixed $h>0$. \\
\\
Thus from \eqref{qq5},\eqref{hder} and \eqref{depl} it follows that,
\begin{align}
\int_{Q_{x_0,h}} d M^{0}_{ij}(x) =\ V^{+}_h(x_0)\hspace{2pt}\underset{\delta\downarrow 0}{lim}\lb\underset{\epsilon \rightarrow 0}{lim\hspace{1pt}sup}\ M^{\theta_{x_0},\epsilon,h-\delta}_{ij}\rb.\label{qq10}
\end{align}
-the last line follows from \eqref{mst} with $(h-\delta)$ in the place of $h$.\\
\\
Now using the fact that $\gamma^{*}_{x_0,\epsilon,h}$ is scale invariant, it follows that
$$ \gamma^{*}_{x_0,\epsilon,h}=\ \gamma^{*}_{x_0,\epsilon,h-\delta}.$$
So, using the relation \eqref{PT and HT relation} it follows that 
$$\theta_{x_0,\epsilon,h}\hspace{1.5pt}M^{\theta_{x_0,\epsilon,h}}=\ \theta_{x_0,\epsilon,h-\delta}\hspace{1.5pt}M^{\theta_{x_0},\epsilon,h-\delta}$$ 
Next using the expression \eqref{qq9} and as it is shown in \eqref{depl}, it simply follows that 
$$ \underset{\delta\downarrow 0}{lim}\hspace{2pt}\underset{\epsilon\rightarrow 0}{lim}\ \frac{\theta_{x_0,\epsilon,h}}{\theta_{x_0,\epsilon,h-\delta}}=\ 1.$$
Thus \eqref{qq10} becomes
$$\int_{Q_{x_0,h}} d M^{0}_{ij}(x) =\ V^{+}_h(x_0)\lb\underset{\epsilon \rightarrow 0}{lim\hspace{1pt}sup}\ M^{\theta_{x_0},\epsilon,h}_{ij}\rb.$$
Since $d M^{0}_{ij}$ is a Radon measure, so it follows that
$$ \int_{\overline{Q}_{x_0,h}} d M^{0}_{ij}(x)\geq\  V^{+}_{h}(x_0)\ \lb\underset{\epsilon \rightarrow 0}{lim\hspace{1pt}sup}\ M^{\theta_{x_0},\epsilon,h}_{ij}\rb. $$
Similarly, one shows ( likewise as we did in Step 2 for $d\mu^{0}$ ) 
$$ \int_{\overline{Q}_{x_0,h}} d M^{0}_{ij}(x)\leq\  V^{-}_{h}(x_0)\ \lb\underset{\epsilon \rightarrow 0}{lim\hspace{1pt}inf}\ M^{\theta_{x_0},\epsilon,h}_{ij}\rb. $$
We consider the sequence of test functions $\psi^\delta_h\in C^{0}(\overline{\Omega})$ defined in \eqref{pshd} and using the fact $d M^{0}$ is a Radon measure we get
$$| \int_{\Omega}\psi^\delta_h d M^{0}_{ij} - \int_{\overline{Q}_{x_0,h}}d M^{0}_{ij}\ |= | \int_{Q_{x_0,h+\delta}\smallsetminus \overline{Q}_{x_0,h}}\psi^\delta_h d M^{0}_{ij}\ | \leq\ |\ d M^{0}_{ij}(Q_{x_0,h+\delta}-\overline{Q}_{x_0,h})\ | \rightarrow 0, \mbox{ as }\delta\rightarrow 0.$$ 
Thus 
\begin{align}
\int_{\overline{Q}_{x_0,h}} d M^{0}_{ij} &=\ \underset{\delta\downarrow 0}{lim}\ \int_{\Omega} \psi^\delta_h d M^{0}_{ij}=\  \underset{\delta\downarrow 0}{lim}\lb \underset{\epsilon\rightarrow 0}{lim}\ \int_{\Omega}\frac{1}{|\omega_{\epsilon}|}\chi_{\omega_{\epsilon}}(x)\nabla w_\epsilon^i(x)\cdot e_j\ \psi^\delta_h dx\rb\notag\\
&=\ \underset{\delta\downarrow 0}{lim}\lb\underset{\epsilon\rightarrow 0}{lim}\ \frac{|\omega_{\epsilon}\cap Q_{x_0,h+\delta}|}{|\omega_{\epsilon}|}\ \frac{1}{|\omega_\epsilon\cap Q_{x_0,h+\delta}|}\int_{Q_{x_0,h+\delta}}\chi_{\omega_{\epsilon}}(x)\nabla \widetilde{w}_{x_0,\epsilon,h+\delta}^i(x)\cdot e_j\ \psi^\delta_h dx\rb.\label{hap}
\end{align}
Again one shows the contribution coming from the annular region $Q_{x_0,h+\delta}\smallsetminus Q_{x_0,h}$ is negligible. We may also replace $|\omega_\epsilon\cap Q_{x_0,h}|$ by 
$|\omega_\epsilon\cap Q_{x_0,h+\delta}|$, similar to as it is shown in \eqref{hder} and \eqref{depl}.
Therefore, from \eqref{hap} it follows that 
\begin{equation}
\int_{\overline{Q}_{x_0,h}} d M^{0}_{ij}(x) =\ V^{-}_h(x_0)\hspace{2pt}\underset{\delta\downarrow 0}{lim\hspace{1.5pt}inf}\lb\underset{\epsilon \rightarrow 0}{lim\hspace{1pt}sup}\ M^{\theta_{x_0},\epsilon,h-\delta}_{ij}\rb. 
\end{equation}
-the last line follows from \eqref{mst} with $h$ in the place of $h+\delta$.\\
Thus \eqref{qq4} follows.
\paragraph{ Step 4 : }
From \eqref{qq4}, it follows, since $V^{+}_h(x_0)=V^{-}_h(x_0)=\mu^{0}(\overline{Q}_{x_0,h})>0$,
$\underset{\epsilon\rightarrow 0}{lim}\ M_{ij}^{\theta_{x_0},\epsilon,h}$ exists and we have
\begin{equation*} \underset{\epsilon\rightarrow 0}{lim}\ M_{ij}^{\theta_{x_0},\epsilon,h} =\ \frac{1}{\mu^{0}(\overline{Q}_{x_0,h})}\int_{\overline{Q}_{x_0,h}} d M^{0}_{ij}(x)\end{equation*}
In this final step, we will be passing to the limit as $h\rightarrow 0$ in the above relation to get the desired result. 
We use the Lebesgue differentiation theorem \cite[Theorem 8.4.6.]{CZ} to have
\begin{equation*}\begin{aligned}
\underset{h\rightarrow 0}{lim}\ \frac{1}{\mu^{0}(\overline{Q}_{x_0,h})}\int_{\overline{Q}_{x_0,h}} d M^{0}_{ij}(x) =\ \underset{h\rightarrow 0}{lim}\ \frac{1}{\mu^{0}(\overline{Q}_{x_0,h})}\int_{\overline{Q}_{x_0,h}} M^{0}_{ij}(x)d\mu^{0}=\ M^{0}_{ij}(x_0), &\\
 \mu^{0}  \mbox{ almost }&\mbox{ everywhere }x_0.
\end{aligned}\end{equation*}
Therefore, we finally get
\begin{equation*} \underset{h \rightarrow 0}{lim}\ \underset{\epsilon \rightarrow 0}{lim}\  M_{ij}^{\theta_{x_0, \epsilon, h}}(x_0)=\ M^0_{ij}(x_0),\quad\mu^{0}\mbox{ almost everywhere }x_0.\end{equation*}
In other words, there exist a sequence $\{\theta_{x_0}^n\}_{n\in\mathbb{N}} \ ( >0 )$ depending upon the point $x_0$,
such that as $n\rightarrow \infty$, $ \theta_{x_0}^n \rightarrow 0 $ and the constant polarization tensors $ M^{\theta_{x_0}^n} \rightarrow M^{0}(x_0)$, $\mu^{0}$ almost everywhere $x_0$. 
\hfill\epr                                                                                                                                                                                   
\begin{remark}
In the above proof, we have dealt with the convergence behavior of the polarization tensor $M^{\theta_{x_0,\epsilon,h}}$. Regarding the associated measure $\mu^{\theta_{x_0,\epsilon,h}}$,
we have 
\begin{equation*}\underset{h\rightarrow 0}{lim}\lb\underset{\epsilon\rightarrow 0}{lim}\ d\mu^{\theta_{x_0,\epsilon,h}}\rb =\ \delta_{x_0},\end{equation*}
because, we have (from \eqref{ED5})
\begin{equation*} d\mu^{\theta_{x_0,\epsilon,h}} =\ \frac{\theta_{x_0,\epsilon,h}}{\int_{Q_{x_0,h}} \theta_{x_0,\epsilon,h}} dx =\ \frac{1}{|Q_{x_0,h}|}dx \quad\mbox{ ( independent of }\epsilon\ ).\end{equation*} 
\hfill\qed
\end{remark}
\noindent
Theorem \ref{approximation} immediately gives the pointwise bounds on the polarization tensor $M^{0}.$ 
\bpr[Proof of Theorem 1.3]
Let's recall from Proposition \ref{bound-mt} the pointwise bounds on the polarization tensor for $M^{\theta}(x)$
for $\mu^{\theta}$ almost everywhere $x\in \Omega$  we have,
\begin{equation}\label{gb-mt}
 \mathrm{min}\left\{ 1, \frac{\gamma_0}{\theta \gamma_0 +(1-\theta)\gamma_1} \right\}|\xi|^2 \leq M^{\theta}_{ij}(x)\xi_i\xi_j \leq \mathrm{max}\left\{ 1, \frac{\gamma_0}{\theta \gamma_0 +(1-\theta)\gamma_1} \right\}|\xi|^2 , \hspace{4pt} \xi \in \mathbb{R}^N. 
\end{equation}
Now take any point $x_0 \in support\ of\ d\mu^{0}$, we want to derive bounds on $M^0(x_0)$, $\mu^{0}$ almost everywhere.
By the previous Theorem \ref{approximation}, there exist a sequence $\theta_{x_0}^n \ ( >0 ) $ depending
upon the point $x_0$ and a sequence of constant polarization tensors $M^{\theta_{x_0}^n}$ such that as 
$n\rightarrow \infty, \theta^n_{x_0}\rightarrow 0$ and $  M^{\theta_{x_0}^n} \rightarrow M^{0}(x_0)$, $\mu^{0}$ almost everywhere $x_0$.
Therefore, we have the following estimate,
\[ \mathrm{min}\left\{ 1, \frac{\gamma_0}{\gamma_1} \right\}|\xi|^2 \leq M^{0}_{ij}(x_0)\xi_i\xi_j \leq \mathrm{max}\left\{ 1, \frac{\gamma_0}{\gamma_1} \right\}|\xi|^2 , \hspace{4pt} \xi \in \mathbb{R}^N. 
\]
\noindent
It immediately shows that $M^{0}(x_0)$ is a positive definite matrix. 
\noindent
Next we recall the optimal bounds on $M^{\theta}(x)$ from the Proposition \ref{Trace bounds - general}. 
$$ \mathrm{trace}\ (\theta^n_{x_0} M^{\theta^n_{x_0}})^{-1}\ \leq\ \frac{N}{\theta^n_{x_0}} - \frac{1-\theta^n_{x_0}}{\theta^n_{x_0}}\Big( 1 - \frac{\gamma_1}{\gamma_0}\Big).$$
This corresponds to the lower curve $c^{\theta^n_{x_0}}_1$ shown in Figure 1. 
Passing to the limit as $n$ tends to infinity, we  simply obtain the lower bound for $M^0(x_0)$ as
$$\mathrm{trace}\ (M^{0}(x_0))^{-1}\ \leq\ (N-1) + \frac{\gamma_1}{\gamma_0}.$$
Similarly for upper bound we recall from Proposition \ref{Trace bounds - general}
$$\ \mathrm{trace}\ \lb I -\theta^n_{x_0}M^{\theta^n_{x_0}}\rb^{-1} \leq \frac{N}{1-\theta^n_{x_0}} + \frac{\theta^n_{x_0}}{1-\theta^n_{x_0}}\left(\frac{\gamma_0}{\gamma_1} - 1\right)$$
This corresponds to the upper curve $c^{\theta^n_{x_0}}_2$ shown in Figure 1.  Since $\theta^n_{x_0} \rightarrow 0$, it follows that
$$\mathrm{trace}\lb I + \theta^n_{x_0} M^{\theta^n_{x_0}}\ \rb\ \leq N(1+ \theta^n_{x_0})\ +\ \theta^n_{x_0}(1-\theta^n_{x_0})\lb \frac{\gamma_0}{\gamma_1} - 1\rb.$$
Thus as $n$ tends to infinity we get the upper bound of $M^0(x_0)$ as
$$\mathrm{trace}\ (M^{0}(x_0)) \leq (N -1) +  \frac{\gamma_0}{\gamma_1}.$$
This ends the proof of Theorem \ref{Optimal trace bounds theorem}. The curves defined by \eqref{LU-bounds} are shown in Figure 1 (in red).  

\hfill\epr
\noindent
In general, for a given density function $\theta\in L^\infty(\Omega;[0,1])$, 
we denote by $\mathcal{G}_\theta$ (the $G$-closure set) the set of all possible H-limits
\begin{align*}
\mathcal{G}_\theta = \{\gamma^{*}\in L^\infty(\Omega;M(\gamma_1,\gamma_0;\Omega))\ |\ \mbox{there exists a characteristic function } \chi_{\omega_\epsilon}(x) \mbox{ satisfying }\eqref{ED7}&\\
\mbox{and }\gamma_\epsilon(x) \ \mbox{ defined by }\eqref{Conductivity_profile}, \mbox{ satisfies }\eqref{ED9}\}.
\end{align*}
Similarly, we define $\mathcal{M}_\theta$ is the set of all possible polarization tensors  corresponding to the density function $\theta(x)$ as
$$\mathcal{M}_\theta = \{M^{\theta} \mbox{ satisfying }\eqref{gb-mt}\ |\ \mbox{ $M^{\theta}$ is defined by \eqref{ED2} and \eqref{ED3}.} \}$$  
And we denote by $\mathcal{M}_0$ the set of all polarization tensors with the near zero volume fraction.
\begin{figure}[H]
 \begin{center}
  \includegraphics[width = 9.8cm]{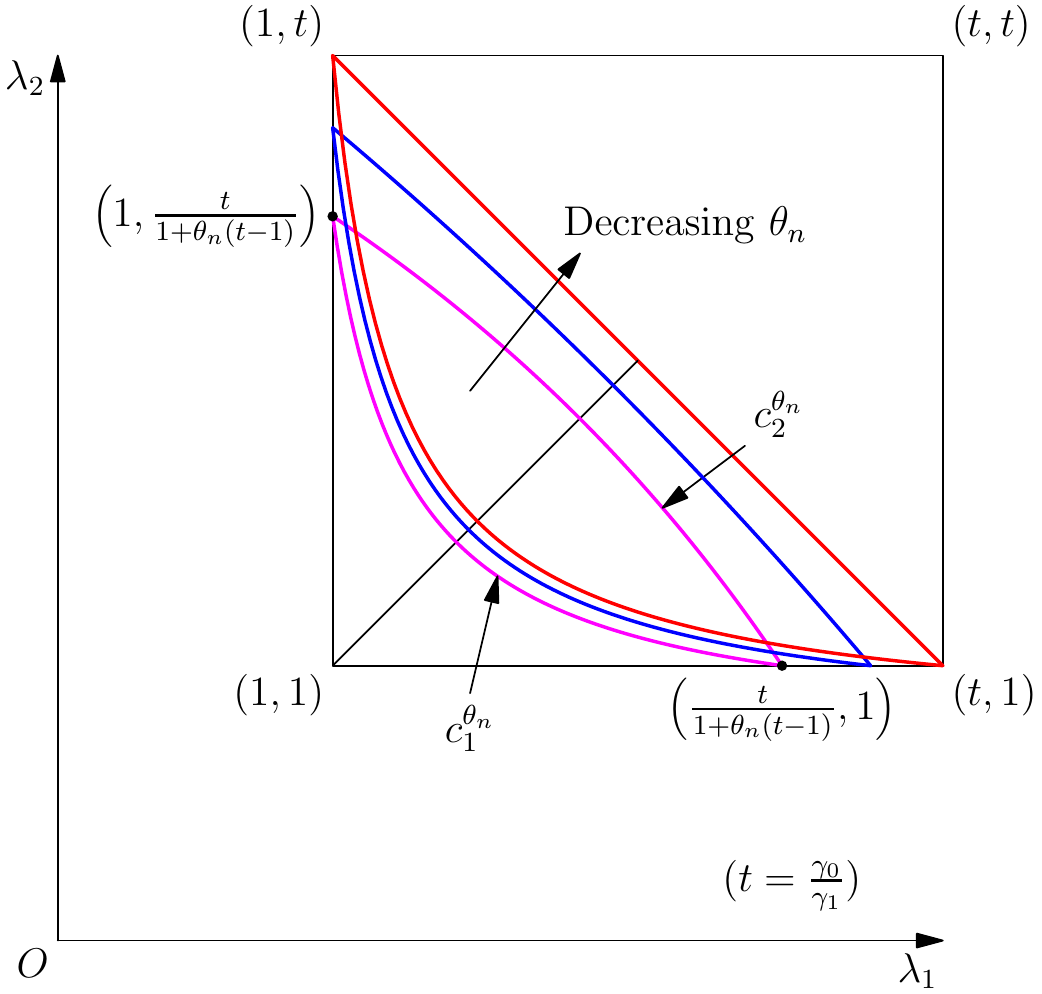}
  \caption{$N=2$: The $\mathcal{M}_0$ region enclosed by pair of red curves can be seen as an approximation of the regions $\mathcal{M}_\theta$ enclosed by the pair of pink and pair of blue curves respectively as $\theta$ goes to zero.}
 \end{center}
\end{figure}
\begin{remark}\label{ED10}
If $x_0$ is such that the polarization tensor
$M^0(x_0)$ satisfies the strict inequality $\textrm{trace}\ (M^{0}(x_0)) < (N-1) + \frac{\gamma_0}{\gamma_1}$,  
then there exists a $\widetilde{\theta}\in (0,1]$ such that $M^{0}(x_0)= M^{\theta}\in \mathcal{M}_{\theta}, \ \forall\ \theta \leq \widetilde{\theta}$.
The existence of such $\widetilde{\theta}$ is clear from the above Figure 1.
\end{remark}
\section{Optimality of the trace bounds for $M^{0}$} \label{Polarization tensors zero volume fraction optimal trace bounds}
\setcounter{equation}{0}
We have derived the bounds for the polarization tensor $M^0$ in Section 3 ( cf. Theorem \ref{Optimal trace bounds theorem}(a) ). 
These bounds were obtained earlier in \cite{CV3,AK}. 
Here we prove the converse namely Theorem \ref{Optimal trace bounds theorem} (b). 
More precisely, we show that given $\{\lambda_1(x),..,\lambda_N(x)\}$ satisfying the
inequalities \eqref{LU-bounds}, there are sequential laminates whose associated
polarization tensor $M^{0}(x)$ has eigen values $\{\lambda_1(x),..,\lambda_N(x)\}$.\\

We begin by computing the polarization tensor for rank-$p$ sequential laminates based on the relationship \eqref{PT and HT relation}. 
The homogenization tensors for such laminates are well-known \cite{A,T}.
\begin{example}[Rank-$p$ Sequential Laminates]\cite[pp. 102]{A}
Let $(e_i)_{1\leq i \leq p}$ be a collection of unit vectors in $\mathbb{R}^{N}$ and $(\theta_i)_{1\leq i \leq p}$ 
the proportions at each stage of the lamination process. Then we have
the following :\\
(a): For a rank-$p$ sequential laminate with matrix $\gamma_0$ and core $\gamma_1$
$$(\prod_{j=1}^{p} \theta_j)( \gamma_p^{*} - \gamma_0I)^{-1}(x) = (\gamma_1-\gamma_0)^{-1}I + \sum_{i=1}^{p}\lb (1-\theta_i)(\prod_{j=1}^{i-1} \theta_j)\rb\frac{e_i\otimes e_i}{\gamma_0}. $$
(b): For rank-$p$ sequential laminate with matrix $\gamma_1$ and core $\gamma_0$
$$(\prod_{j=1}^{p} (1-\theta_j))( \gamma_p^{*} - \gamma_1I)^{-1}(x) = (\gamma_1-\gamma_0)^{-1}I + \sum_{i=1}^{p}\lb \theta_i(\prod_{j=1}^{i-1}(1-\theta_j)\rb\frac{e_i\otimes e_i}{\gamma_1}.$$
The following lemma from is important for proving the optimality. 
\end{example}
\begin{lemma}\cite{A}\label{lam lemma}
Let  $(e_i)_{1\leq i \leq p}$ be a collection of unit vectors. Let $\theta\in (0,1]$. Now for a fixed point $x$ and for any collection of non-negative real numbers  $(m_i)_{1\leq i \leq p}$ satisfying $\sum_{i=1}^{p} m_i =1 $, 
there exists a rank-p sequential laminate $\gamma^{*}_p(x)$ with matrix $\gamma_0I$ and core $\gamma_1I$ in proportion $(1-\theta)$ and $\theta$ respectively and with lamination directions $(e_i)_{1\leq i \leq p}$ such that
\begin{equation}\label{Matrix 0, Core 1}
\theta( \gamma_p^{*} - \gamma_0I)^{-1}(x) = (\gamma_1-\gamma_0)^{-1}I + (1-\theta)\sum_{i=1}^{p} m_i\frac{e_i\otimes e_i}{\gamma_0e_i.e_i}.
\end{equation}
\end{lemma}
An analogous result holds when the roles of $\gamma_0$ and $\gamma_1$ (in proportions $(1-\theta)$ and $\theta$ respectively) in the lemma above are switched. The formula above is replaced by
\begin{equation}\label{Matrix 1, Core 0}
(1-\theta)(\gamma_p^{*} - \gamma_1I)^{-1}(x) = (\gamma_0-\gamma_1)^{-1}I + \theta\sum_{i=1}^{p} m_i\frac{e_i\otimes e_i}{\gamma_1e_i.e_i}.
\end{equation}
\hfill\qed\\
Denoting the polarization tensor $M^{\theta}_{p}$ in the case of matrix $\gamma_0I$ and core $\gamma_1I$, by (\ref{PT and HT relation}) and \eqref{Matrix 0, Core 1}  we then have,
\begin{equation}\label{lam-p1}\lb M_p^{\theta}(x)\rb^{-1} = I + (1-\theta)\frac{(\gamma_1-\gamma_0)}{\gamma_0}\sum_{i=1}^{p} m_i (e_i\otimes e_i).\end{equation}

For matrix $\gamma_1I$ and core $\gamma_0I$ from \eqref{PT and HT relation} and \eqref{Matrix 1, Core 0}, we have
\begin{equation}\label{lam-p2}
(I-\theta M^{\theta}_p)(x)^{-1}= \frac{1}{(1-\theta)} I + \frac{\theta}{1-\theta}\frac{(\gamma_0-\gamma_1)}{\gamma_1}\sum_{i=1}^p m_i (e_i\otimes e_i)
\end{equation}
\hfill\qed\\
\\
\bpr[Proof of Theorem \ref{Optimal trace bounds theorem}]
Let $M(x_0)$ be a second order positive definite tensor with its eigenvalues $(\lambda_1(x_0),..,\lambda_N(x_0))$ lying on the curve defining the upper bound in the Figure 1 :  
\begin{equation}\label{equality} 1\leq \lambda_i(x_0) \leq \frac{\gamma_0}{\gamma_1}, \quad1\leq i \leq N \quad \mbox{ and }\quad \sum_{i=1}^N \lambda_i(x_0) = (N-1) + \frac{\gamma_0}{\gamma_1}.\end{equation}
Then there exists a collection of non-negative real numbers  $(m_i)_{1\leq i \leq N}$ depending upon the point $x_0$ satisfying 
\begin{equation}\label{m_i}
\sum_{i=1}^{N} m_i =1\quad \mbox{ and }\quad \lambda_i(x_0) = 1 + m_i(\frac{\gamma_0}{\gamma_1}-1 ).
\end{equation}
After getting such $(m_i)_{1\leq i \leq N}$\ 's through the equation \eqref{m_i}, and choosing $\theta\in (0,1]$
we consider the  rank-$N$ sequential laminated structure defined in the Lemma \ref{lam lemma}.
We obtain a polarization tensor in  $M^{\theta}_N $ with 
matrix $\gamma_1$ and core $\gamma_0$ such that
$$(I-\theta M^{\theta}_N)(x_0)^{-1} = \frac{1}{(1-\theta)} I + \frac{\theta}{1-\theta}\frac{(\gamma_0-\gamma_1)}{\gamma_1}\sum_{i=1}^N m_i( e_i\otimes e_i).$$
The eigenvalues of $M^{\theta}_N(x_0)$ denoted by $(\lambda^{\theta}_i(x_0))_{1\leq i \leq N}$
are given by
\begin{equation}\label{lam m}
\frac{1}{1-\theta\lambda^{\theta}_i(x_0)} = \frac{1}{(1-\theta)} + \frac{\theta}{1-\theta}\frac{(\gamma_0-\gamma_1)}{\gamma_1} m_i.
\end{equation}
Since $m_i$ satisfies $0\leq m_i\leq1$, it follows
$$ 1\ \leq\ \lambda_i^{\theta}(x_0)\ \leq\ \frac{\gamma_0}{\theta (\gamma_0-\gamma_1) + \gamma_1} $$
and since $\sum_{i=1}^N m_i =1$, we have also 
$$\sum_{i=1}^N (1-\theta\lambda_i^{\theta}(x_0))^{-1} = \frac{N}{(1-\theta)} + \frac{\theta}{1-\theta}\frac{(\gamma_0-\gamma_1)}{\gamma_1}.$$

Now choosing a sequence $\theta_n\rightarrow 0$ ( e.g: $\theta_n= \frac{1}{n}),$ we will show as $n\rightarrow \infty$ $\lambda_i^{\theta_n}(x_0)\rightarrow \lambda_{i}(x_0).$
From the asymptotic expansion of the equation \eqref{lam m} in terms of $\theta\approx 0$ we have
$$(1 + \theta_n \lambda_i^{\theta_n}(x_0))\ \approx\ (1+\theta_n) + \theta_n(1+\theta_n)\frac{(\gamma_0-\gamma_1)}{\gamma_1} m_i.$$
So as $n\rightarrow \infty$,
$$\lambda^{\theta_n}_i(x_0) \rightarrow\  1 + \frac{(\gamma_0-\gamma_1)}{\gamma_1}m_i =\lambda_i(x_0).$$
This is equivalent to saying 
$$M^{\theta_n}_N(x_0) \rightarrow M(x_0).$$
To finish the proof we need to show $M\in \mathcal{M}_0$ ( see the end of Section 3 for the definition of $\mathcal{M}_0$ ) 
To this end, since we are working at the point $x_0$, we can in fact assume that we are dealing with periodic rank- $N$ sequntial laminates.
Now the result quickly follows \cite[Theorem 8.1]{AK}.

Similarly one shows that equality of lower bound can be achieved through the $N$ sequential laminates with matrix $\gamma_0$ and core $\gamma_1$.

To show that any interior point in the region defined by the bounds \eqref{LU-bounds} corresponds to a polarization tensor of near zero volume fraction, we can 
follow the arguments found in \cite[Page no.124-125]{A}.

Thus, any tensor $M(x)$ satisfying the pointwise bounds given by \eqref{LU-bounds} is in $\mathcal{M}_{0}$, which completes our discussion on optimality of the bounds. 
\hfill\epr

\paragraph{Acknowledgement :}
This work has been carried out within a project supported by the 
Airbus Group Corporate Foundation Chair ``Mathematics of Complex Systems'' 
established at Tata Institute Of fundamental Research (TIFR) - Centre for Applicable Mathematics.

\bibliographystyle{plain}
\bibliography{Master_bibfile}

\end{document}